\documentclass[12pt,reqno]{amsart}
\usepackage{mathrsfs}
\usepackage{graphicx,epsfig,amsmath,amsfonts,dsfont}
\usepackage{color}
\usepackage{latexsym,amssymb}
\usepackage{lineno,version}
\usepackage{multirow}
\usepackage{bm}
\usepackage{float}
\usepackage{booktabs}
\usepackage{hyperref} % link
\usepackage{geometry}
\allowdisplaybreaks

\catcode`\@=11
\theoremstyle{plain}
\@addtoreset{equation}{section}   % Makes \section reset 'equation' counter.

\@addtoreset{figure}{section}
\renewcommand\thefigure{\thesection.\@arabic\c@figure}
\newtheorem{thm}{\noindent Theorem}[section]
\newtheorem{rem}{\noindent Remark}[section]

\newtheorem{lem}{\noindent Lemma}[section]
\newtheorem{examp}{\noindent Example}[section]

\newcommand{\ba}{\begin{array}}\newcommand{\ea}{\end{array}}
\newcommand{\be}{\begin{eqnarray}}\newcommand{\ee}{\end{eqnarray}}
\newcommand{\beq}{\begin{equation}}\newcommand{\eeq}{\end{equation}}
\newcommand{\bex}{\begin{eqnarray*}}
\newcommand{\eex}{\end{eqnarray*}}

\font\tenbi=cmmib10   at 11 pt
\font\sevenbi=cmmib10 at 9pt
\font\fivebi=cmmib7 at 6pt
\newfam\bifam
\textfont\bifam=\tenbi \scriptfont\bifam=\sevenbi
\scriptscriptfont\bifam=\fivebi

\font\tendb=msbm10 at 12 pt
\font\sevendb=msbm7
\newfam\dbfam
\textfont\dbfam=\tendb \scriptfont\dbfam=\sevendb

\begin{document}

\title[Efficient numerical method for multi-term FDEs with C-F derivatives]
{Efficient numerical method for multi-term time-fractional diffusion equations with Caputo-Fabrizio
derivatives$^*$}
\author[]
{Bin Fan$^{1,\dag}$}
\thanks{\hskip -12pt
${}^*$This work was supported by the Natural Science Foundation of Fujian Province (Grant No. 2022J05198), the Foundation of Fujian Provincial Department of Education (Grant No. JAT210293), and the Natural Science Foundation of Fujian University of Technology (Grant No. GY-Z21069).\\
${}^{1}$School of Computer Science and Mathematics, Fujian University of Technology, 350118 Fuzhou, China.\\
Email: bfan@fjut.edu.cn}

%\subjclass[2010]

%\date{\today}

\maketitle

\begin{abstract}
In this paper, we consider a numerical method for the multi-term Caputo-Fabrizio time-fractional diffusion equations (with orders $\alpha_i\in(0,1)$, $i=1,2,\cdots,n$).
The proposed method employs a fast finite difference scheme to approximate multi-term fractional derivatives in time, requiring only $O(1)$ storage and $O(N_T)$ computational complexity, where $N_T$ denotes the total number of time steps. Then we use a Legendre
spectral collocation method for spatial discretization. The stability and convergence of the scheme have been thoroughly discussed and rigorously established.
We demonstrate that the proposed scheme is unconditionally stable and convergent with an order of $O\left(\left(\Delta t\right)^{2}+N^{-m}\right)$, where $\Delta t$, $N$, and $m$ represent the timestep size, polynomial degree, and regularity in the spatial variable of the exact solution, respectively.
Numerical results are presented to validate the theoretical predictions.

\end{abstract}

\section{Introduction}\label{sec:1}

Fractional differential equations  have wide applications in various fields of science, including physics, economics, engineering, chemistry, biology and others\cite{cuesta2012image,magin2004fractional,mainardi2022fractional,metzler2000random,zhou2023BS}. There are many kinds of definitions for the fractional derivatives,
the most used fractional derivatives are the Riemann-Liouville fractional derivative and the Caputo fractional derivative \cite{podlubny1998fractional,jin2021fractional}. However, both of these operators still present challenges in practical applications.
To be more precise, the Riemann-Liouville derivative of a constant is non-zero and the Laplace transform of this derivative contains terms that lack physical significance. The Caputo fractional derivative has successfully addressed both issues, however, its definition involves a singular kernel which poses challenges in analysis and computation. Caputo and Fabrizio \cite{caputo2015new} have proposed a novel definition of the fractional derivative with a smooth kernel, referred to as the Caputo and Fabrizio (CF) derivatives, which present distinct representations for the temporal and spatial variables. The representation in time variable
is suitable to use the Laplace transform, and the spatial representation is more convenient to use the Fourier transform.
Although there is ongoing debate regarding the mathematical properties of fractional derivatives with non-singular kernels \cite{2020Why,2020Fractional}, numerous scholars remain interested in studying differential equations involving such derivatives due to their nice performance in various applications.
Considering the CF derivative offers two primary advantages: 1) The utilization of a regular kernel in non-local systems is motivated by its potential to accurately depict material heterogeneities and fluctuations of various scales, which cannot be adequately captured by classical local theories or fractional models with singular kernels, see, e.g., \cite{caputo2015new,caputo2016applications};
2) CF derivatives have numerical advantages. As we know, the truncation error of the numerical calculation for  fractional operators with singular kernels is typically dependent on the order $\alpha$. For instance, in the case of the Caputo fractional derivative, employing classical L1 discretization results in an error of order $2-\alpha$, which becomes highly unfavorable when $\alpha\approx1$. In order to enhance actuarial accuracy, the utilization of higher-order methods will lead to an increase in computational complexity, particularly for problems with high dimensions. However, within the same approximation framework, the CF derivative has a higher truncation error, see Remark \ref{rem-trun-err}.
Further properties and diverse applications of this fractional derivative can be found in various references, such as \cite{caputo2016applications,gomez2016modeling,atangana2016new,atangana2016new-2,mirza2017fundamental,tuan2020well,jia2023analysis}.

Let $T,L>0$, $\Lambda:=(0,S)$. In this paper, we are concerned with the numerical approximation of the multi-term time-fractional diffusion equation
\begin{equation}\label{pro}
  P_{\alpha_1,\alpha_2,\cdots,\alpha_{n-1},\alpha_n}\left(_0^{\mathrm{CF}}\!D_t\right)u(x,t)=\frac{\partial^2u(x,t)}{\partial x^2}+f(x,t),\quad (x,t)\in\Lambda\times(0,T],
\end{equation}
with the initial conditions
\begin{equation}\label{pro-init-cond}
 u(x,0)=\varphi(x),\quad x\in\Lambda,
\end{equation}
and the boundary condition
\begin{equation}\label{pro-bou-cond}
  u(0,t)=u(S,t)=0, \quad t\in[0,T],
\end{equation}
where
\begin{equation}\label{multi-term-cap}
  P_{\alpha_1,\alpha_2,\cdots,\alpha_{n-1},\alpha_n}\left(_0^{\mathrm{CF}}\!D_t\right)u(x,t)=\sum\limits_{i=1}^nd_i\cdot {}_0^{\mathrm{CF}}\!D_t^{\alpha_{i}}u(x,t),
\end{equation}
$0<\alpha_1\leq\alpha_2\leq\cdots\leq\alpha_{n-1}\leq\alpha_n<1$, and $d_i\geq0,~i=1,2,\cdots,n, n\in\mathds{N}$. $\varphi(x)$ and $f(x,t)$ are given sufficiently smooth functions in their respective domains.
In addition, ${}_0^{\mathrm{CF}}\!D_t^{\alpha_i}u(x,t)$ is the Caputo-Fabrizio derivative
   operator of order $\alpha_i$ \cite{caputo2015new,caputo2016applications} defined as
\begin{eqnarray}
   &&  {}_0^{\mathrm{CF}}\!D_t^{\alpha}u(x,t)=\frac{1}{1-\alpha}\int_0^{t}\frac{\partial u(x,s)}{\partial s}\exp\left(-\alpha\frac{t-s}{1-\alpha}\right)\mathrm{d}s. \label{def-pCFD}
\end{eqnarray}
If $n=1$, then $(\ref{pro})$-$(\ref{pro-bou-cond})$ reduces to the single-term time-fractional diffusion equation.
The model of $(\ref{pro})$-$(\ref{pro-bou-cond})$, which describes the temporal flow of water within a leaky aquifer at various scales \cite{atangana2016new-2,Djida2017More}, as well as the electro-magneto-hydrodynamic flow of non-Newtonian biofluids with heat transfer \cite{2017Modeling}, etc. For the well-posedness of $(\ref{pro})$-$(\ref{pro-bou-cond})$, we refer to, e.g., \cite{Mohammed2017,Nasser2016,tuan2020well}.

Many researchers have explored the numerical approximation of both single-term and multi-term time fractional diffusion equations.
In \cite{liu2004analysis}, Liu et al. proposed a finite difference method for solving time-fractional diffusion equations in both space and time domains.
Lin and Xu \cite{lin2007finite} utilized a finite difference scheme in time and Legendre spectral methods in space to numerically solve the time-fractional diffusion equations. Subsequently, Li and Xu \cite{li2009space}improved upon their previous work by proposing a space-time spectral method for these equations.
For the numerical treatment of multi-term time-fractional diffusion equations, \cite{ren2014efficient} proposed a fully-discrete schemes for one- and two-dimensional multi-term time fractional sub-diffusion equations. These schemes combine the compact difference method for spatial discretization with L1 approximation for time discretization.
The Galerkin finite element method and the spectral method were introduced in \cite{jin2015galerkin} and \cite{zheng2016high,fardi2022legendre}, respectively.
Zhao et al. \cite{zhao2017convergence} developed a fully-discrete scheme for a class of two-dimensional multi-term time-fractional diffusion equations with Caputo fractional derivatives, utilizing the finite element method in spatial direction and classical L1 approximation in temporal direction.
Akman et al. \cite{akman2018new} proposed a numerical approximation called the L1-2 formula for the Caputo-Fabrizio derivative using quadratic interpolation.
In \cite{Yu2019}, finite difference/spectral approximations for solving two-dimensional time CF fractional diffusion equation were proposed and analyzed. Later, a second order scheme \cite{shi2020} was devised for addressing this problem.
A compact alternating direction implicit (ADI) difference scheme was proposed by \cite{taghipour2020new} for solving the two-dimensional time-fractional diffusion equation.

Simulating models with fractional derivatives presents a challenge due to their non-locality, which significantly impedes algorithm efficiency and necessitates greater memory storage compared to traditional local models.
In particular, for fractional models, the computational complexity of obtaining an approximate solution is $O(N_T^2)$ and the required memory storage is $O(N_T)$, which contrasts with local models that have a complexity of $O(N_T)$ and require a memory storage of $O(1)$, where $N_T$ denotes the total number of time steps, see, e.g., \cite{lin2007finite,Yu2019,shi2020}.
To address this issue, several researchers have proposed efficient algorithms for computing the derivatives of Riemann-Liouville, Caputo, and Riesz fractional operators, see e.g., \cite{jiang2017fast,zeng2018stable,li2018fast,zhu2019fast} and the references therein.
Recently, a fast compact finite difference method for quasi-linear time-fractional parabolic equations is presented and analyzed in \cite{liu2019fast}. Then, \cite{liu2020fast} proposed a fast second-order numerical scheme for approximating the Caputo-Fabrizio fractional derivative at node $t=t_{k+1/2}$ with computational complexity of $O(N_T)$ and memory storage of $O(1)$.

Inspired by the above mentioned, we extend finite difference/spectral approximations for the multi-term Caputo-Fabrizio time-fractional diffusion equation \eqref{pro}-\eqref{pro-bou-cond}.
Firstly, we present a L1 formula for the Caputo-Fabrizio derivative. In this context, we introduce two discrete fractional differential operators, namely $L_t^{\alpha}$ and $F_t^{\alpha}$, which are essentially equivalent.
However, $F_t^{\alpha}$ effectively utilizes the exponential kernel and incurs lower storage and computational costs compared to $L_t^{\alpha}$.
The idea of this approach is essentially identical to that of reference \cite{liu2020fast}, albeit with a slightly different formulation in our case; specifically, the approximation is centered at point $t=t_k$ and presented in a more concise manner. The error bounds associated with these two operators will be examined in detail.
Secondly, we develop a semi-discrete scheme based on finite difference method for multi-term time-fractional derivatives, with complete proofs of its unconditional stability and convergence rate.
A detailed error analysis is carried out for the semi-discrete problem, showing that the temporal accuracy is second order.
Finally, we present the fully-discrete scheme based on the Legendre spectral collocation method for spatial discretization. We will investigate both the convergence order of this method and its implementation efficiency, while providing a rigorous proof of its spectral convergence in this paper.

The rest of this paper is organized as follows. In Section 2, a semi-discrete scheme is proposed for \eqref{pro}-\eqref{pro-bou-cond} based on fast L1 finite difference scheme.  The stability and convergence analysis of the semi-discrete scheme is presented. In Section 3,
we construct a Legendre spectral collocation method for the spatial discretization of the semi-discrete scheme.
Error estimates are provided for the full discrete problem. Some numerical results are reported in Section 4. Finally, the conclusions are given in Section 5.

\section{Semi-discretization}
Define $t_k:=k\Delta t$, $k=0,1,\cdots,N_T$, where $\Delta t:=T/N_T$ is the time step.
\subsection{Fast L1 formula for Caputo-Fabrizio derivative}\label{sec:2}
We first give L1 approximation for fractional Caputo-Fabrizio derivative of function $h(t)$ defined by
\begin{eqnarray}
   &&  {}_0^{\mathrm{CF}}\!D_t^{\alpha}h(t)=\frac{1}{1-\alpha}\int_0^{t}h'(s)\exp\left(-\alpha\frac{t_k-s}{1-\alpha}\right)\mathrm{d}s.\label{def-CFD}
\end{eqnarray}
In order to simplify the notations, we denote $h(t_k):=h^k$ for $0\leq k\leq M_t$. The L1 formula is obtained by substituting the linear Lagrange interpolation of $h(t)$ into (\ref{def-CFD}). Precisely, the linear approximation of the function $h(t)$ on $[t_{j-1},t_j]$ is written as
\begin{eqnarray}
   && \Pi_{1,j}h(t)=\frac{t_j-t}{\Delta t}h^{j-1}+\frac{t-t_{j-1}}{\Delta t}h^j,\quad 1\leq j\leq k,\label{line-inter}
\end{eqnarray}
and the error in the approximation is
\begin{eqnarray}
   && h(t)-\Pi_{1,j}h(t)=\frac{1}{2}h''(\xi_j)(t-t_{j-1})(t-t_j),\quad  \xi_j\in (t_{j-1},t_j),\quad 1\leq j\leq k.\label{err-line-inter}
\end{eqnarray}
Then we define the discrete fractional differential operator $L_t^\alpha$ by
\begin{align}
    L_t^\alpha h^k&:=\frac{1}{1-\alpha}\sum_{j=1}^k\int_{t_{j-1}}^{t_j}\left(\Pi_{1,j}h(s)\right)'\exp\left(-\alpha\frac{t_k-s}{1-\alpha}\right)\mathrm{d}s \nonumber \\
&=\frac{1}{1-\alpha}\sum_{j=1}^k\frac{h^j-h^{j-1}}{\Delta t}\int_{t_{j-1}}^{t_j}\exp\left(-\alpha\frac{t_k-s}{1-\alpha}\right)\mathrm{d}s \nonumber \\
&=\frac{1}{\alpha\Delta t}\sum_{j=1}^k\left(h^j-h^{j-1}\right)\left(\sigma_{j,k}-\sigma_{j-1,k}\right)\nonumber \\
&=\frac{1}{\alpha\Delta t}\sum_{j=1}^kb_{j,k}\left(h^j-h^{j-1}\right)\nonumber \\
&=\frac{1}{\alpha\Delta t}\left(b_{k,k}h^k+\sum_{j=1}^{k-1}\left(b_{j,k}-b_{j+1,k}\right)h^j-b_{1,k}h^0\right),\label{L1-1-uk}
\end{align}
where $b_{j,k}:=\sigma_{j,k}-\sigma_{j-1,k}$ and
\begin{eqnarray*}
   && \sigma_{j,k}:=\exp\left(-\alpha\frac{t_{k}-t_j}{1-\alpha}\right),\quad 1\leq k\leq N_T,\quad 1\leq j\leq k.
\end{eqnarray*}

The right hand side of \eqref{L1-1-uk} involves a sum of all previous solutions $\left\{h^j\right\}_{j=0}^k$, which reflects the memory effect of the non-local fractional derivative. Thus it requires on average $O\left(N_T\right)$ storage and the total computational cost is $O\left(N_T^2\right)$ with $N_T$ the total number of time steps. This makes both the computation and memory expensive, specially in the case of
long time integration. In order to overcome this difficulty, we propose a further approach to the fractional derivative.
The idea consists in first splitting the convolution integral in \eqref{def-CFD} into a sum of history part and local part as follows:
\begin{align}
     {}_0^{\mathrm{CF}}\!D_t^{\alpha}h^k&=\frac{1}{1-\alpha}\int_0^{t_k}h'(s)\exp\left(-\alpha\frac{t_k-s}{1-\alpha}\right)\mathrm{d}s \nonumber\\
  &=\frac{1}{1-\alpha}\int_0^{t_{k-1}}h'(s)\exp\left(-\alpha\frac{t_k-s}{1-\alpha}\right)\mathrm{d}s+\frac{1}{1-\alpha}\int_{t_{k-1}}^{t_{k}}h'(s)\exp\left(-\alpha\frac{t_k-s}{1-\alpha}\right)\mathrm{d}s  \nonumber\\
  &:=C_h(t_k)+C_l(t_k). \nonumber
\end{align}
Note that a comparable treatment is employed in reference \cite{jiang2017fast}.
Then the history part $C_h(t_k)$ can be rewritten as
\begin{align}
    C_h(t_k)&=\frac{1}{1-\alpha}\int_0^{t_{k-1}}h'(s)\exp\left(-\alpha\frac{t_{k-1}-s+t_k-t_{k-1}}{1-\alpha}\right)\mathrm{d}s \nonumber\\
   &=\exp\left(-\frac{\alpha\Delta t}{1-\alpha}\right)\frac{1}{1-\alpha}\int_0^{t_{k-1}}h'(s)\exp\left(-\alpha\frac{t_{k-1}-s}{1-\alpha}\right)\mathrm{d}s \nonumber\\
   &=\exp\left(-\frac{\alpha\Delta t}{1-\alpha}\right) {}_0^{\mathrm{CF}}\!D_t^{\alpha}h^{k-1},\nonumber
\end{align}
hence we have
\begin{eqnarray}
   && {}_0^{\mathrm{CF}}\!D_t^{\alpha}h^k=\exp\left(-\frac{\alpha\Delta t}{1-\alpha}\right) {}_0^{\mathrm{CF}}\!D_t^{\alpha}h^{k-1}+C_l(t_k). \label{rec-CFD}
\end{eqnarray}
Using the simple recurrence relation \eqref{rec-CFD}, we define the discrete fractional differential operator $F_t^{\alpha}$ by
\begin{align}
   F_t^{\alpha}h^1&=\frac{1}{1-\alpha}\frac{h^1-h^0}{\Delta t}\int_0^{t_{1}}\exp\left(-\alpha\frac{t_{1}-s}{1-\alpha}\right)\mathrm{d}s \nonumber \\
    &=\frac{1}{\alpha \Delta t}\left(\sigma_{1,1}-\sigma_{0,1}\right)\left(h^1-h^0\right)=\frac{b_{1,1}}{\alpha \Delta t}\left(h^1-h^0\right),\label{app-CFD-rec-t1}\\[6pt]
    F_t^{\alpha}h^k&=\exp\left(-\frac{\alpha\Delta t}{1-\alpha}\right)F_t^{\alpha}h^{k-1}+\frac{1}{1-\alpha}\frac{h^k-h^{k-1}}{\Delta t}\int_{t_{k-1}}^{t_{k}}\exp\left(-\alpha\frac{t_{k}-s}{1-\alpha}\right)\mathrm{d}s\nonumber \\
    &=\exp\left(-\frac{\alpha\Delta t}{1-\alpha}\right)F_t^{\alpha}h^{k-1}+\frac{1}{\alpha \Delta t}\left(\sigma_{k,k}-\sigma_{k-1,k}\right)\left(h^k-h^{k-1}\right)\nonumber \\
    &=\exp\left(-\frac{\alpha\Delta t}{1-\alpha}\right)F_t^{\alpha}h^{k-1}+\frac{b_{k,k}}{\alpha \Delta t}\left(h^k-h^{k-1}\right),\qquad k\geq 2.\label{app-CFD-rec-tk}
\end{align}
It is not difficult to see that $F_t^{\alpha}h^k=L_t^{\alpha}h^k$ for $1\leq k\leq N_T$. Comparing $L_t^\alpha h^k$ in (\ref{L1-1-uk}) with $F_t^\alpha h^k$ in \eqref{app-CFD-rec-t1}-\eqref{app-CFD-rec-tk}, the former requires all the previous time step values of $h(t)$ while the latter only needs $h^{k-1}$, $h^{k}$ and $F_t^{\alpha}h^{k-1}$. This implies that approximating ${}_0^{\mathrm{CF}}\!D_t^{\alpha}h^k$ by $F_t^\alpha h^k$ considerably
reduces the storage and computational costs as compared to $L_t^\alpha h^k$. Roughly speaking, replacing $L_t^\alpha h^k$ by $F_t^\alpha h^k$ allows to reduce the storage cost
from $O\left(N_T\right)$ to $O\left(1\right)$, and the computational cost from $O\left(N_T^2\right)$ to $O\left(N_T\right)$.
\begin{rem}
The fast algorithm of Caputo derivative in \cite{jiang2017fast} should be noted to retain an additional truncation error $\varepsilon$, whereas the fast algorithm of CF derivative does not introduce this error. Furthermore, it is worth mentioning that other algorithms, such as parallel computational methods \cite{Gu2020}, result in an augmented spatial complexity.
\end{rem}

The following lemma provides an error bound for approximation $F_t^\alpha h^k$.
\begin{lem}\label{lem-erro-Ft}
Suppose that $h(t)\in \mathds{C}^2[0,T]$. For any $0<\alpha<1$, let
\begin{eqnarray}
   && R_k:={}_0^{\mathrm{CF}}\!D_t^{\alpha}h^k-F_t^\alpha h^k. \nonumber
\end{eqnarray}
Then
\begin{eqnarray}
   && \left|R_k\right|\leq\frac{\alpha T\max\limits_{t\in[0,T]}h''(t)}{8(1-\alpha)^2}\left(\Delta t\right)^2,\qquad j=1,2,\cdots M_t.\nonumber
\end{eqnarray}
\end{lem}
\noindent\emph{Proof.} We consider proving the following estimate by mathematical induction:
\begin{eqnarray}
   &&  \left|R_j\right|\leq\frac{\alpha\max\limits_{t\in[0,T]}h''(t)}{8(1-\alpha)^2}j\left(\Delta t\right)^3,\qquad j=1,2,\cdots M_t.\label{lem-erro-Ft-eq1}
\end{eqnarray}
First we have
\begin{align}
     \left|R_1\right|&=\left|{}_0^{\mathrm{CF}}\!D_t^{\alpha}h^1-F_t^\alpha h^1\right|=\left|\frac{1}{1-\alpha}\int_0^{t_1}\left(h(s)-\Pi_{1,1}h(s)\right)'\exp\left(-\alpha\frac{t_1-s}{1-\alpha}\right)\mathrm{d}s\right| \nonumber\\
    &\leq\left|\frac{1}{1-\alpha}\left(h(s)-\Pi_{1,1}h(s)\right)\exp\left(-\alpha\frac{t_1-s}{1-\alpha}\right)\bigg|_{s=0}^{s=t_1}\right|\nonumber\\
    &\qquad+\left|\frac{\alpha}{(1-\alpha)^2}\int_0^{t_1}\left(h(s)-\Pi_{1,1}h(s)\right)\exp\left(-\alpha\frac{t_1-s}{1-\alpha}\right)\mathrm{d}s\right| \qquad\text{(Integration by parts)}\nonumber\\
    &=\left|\frac{1}{1-\alpha}\frac{1}{2}h''(\xi_1)(s-t_{0})(s-t_1)\exp\left(-\alpha\frac{t_1-s}{1-\alpha}\right)\bigg|_{s=0}^{s=t_1}\right|\nonumber\\
    &\qquad+\left|\frac{\alpha }{(1-\alpha)^2}\int_0^{t_1}\frac{1}{2}h''(\xi_1)(s-t_{0})(s-t_1)\exp\left(-\alpha\frac{t_1-s}{1-\alpha}\right)\mathrm{d}s\right| \qquad (\text{By \eqref{err-line-inter}})\nonumber\\
  &=\left|\frac{\alpha }{2(1-\alpha)^2}\int_0^{t_1}h''(\xi_1)(s-t_{0})(s-t_1)\exp\left(-\alpha\frac{t_1-s}{1-\alpha}\right)\mathrm{d}s\right|\nonumber\\
   &\leq\frac{\alpha}{2(1-\alpha)^2}\max_{t\in[0,T]}h''(t)\frac{1}{4}\left(\Delta t\right)^2\int_0^{t_1}\exp\left(-\alpha\frac{t_1-s}{1-\alpha}\right)\mathrm{d}s\nonumber\\
    &=\frac{\alpha}{2(1-\alpha)^2} \max_{t\in[0,T]}h''(t) \frac{1}{4}\left(\Delta t\right)^2\exp\left(-\alpha\frac{t_1-\eta_1}{1-\alpha}\right)\Delta t \qquad\text{(Mean value theorem)}\nonumber\\[6pt]
    &\leq\frac{\alpha\max\limits_{t\in[0,T]}h''(t)}{8(1-\alpha)^2}\left(\Delta t\right)^3,\nonumber
\end{align}
where $\eta_1\in(0,t_1)$. Therefore, (\ref{lem-erro-Ft-eq1}) holds for $j=1$. Now suppose that (\ref{lem-erro-Ft-eq1}) holds for $j=k-1$, we need to prove that it holds also for $j=k$. Similar to the proof of $\left|R_1\right|$, we can easily get that
\begin{eqnarray}
   && \left|\frac{1}{1-\alpha}\int_{t_{k-1}}^{t_k}\left(h(s)-\Pi_{1,k}h(s)\right)'\exp\left(-\alpha\frac{t_k-s}{1-\alpha}\right)\mathrm{d}s\right|\leq\frac{\alpha\max\limits_{t\in[0,T]}h''(t)}{8(1-\alpha)^2}\left(\Delta t\right)^3.\nonumber
\end{eqnarray}
By combining (\ref{rec-CFD}) and (\ref{app-CFD-rec-tk}), we obtain
\begin{align}
     \left|R_{k}\right|&\leq\exp\left(-\frac{\alpha\Delta t}{1-\alpha}\right)\left|{}_0^{\mathrm{CF}}\!D_t^{\alpha}h^{k-1}-F_t^{\alpha}h^{k-1}\right|+\left|\frac{1}{1-\alpha}\int_{t_{k-1}}^{t_k}\left(h(s)-\Pi_{1,k}h(s)\right)'\exp\left(-\alpha\frac{t_k-s}{1-\alpha}\right)\mathrm{d}s\right| \nonumber\\
   &\leq\left|{}_0^{\mathrm{CF}}\!D_t^{\alpha}h^{k-1}-F_t^{\alpha}h^{k-1}\right|+\left|\frac{1}{1-\alpha}\int_{t_{k-1}}^{t_k}\left(h(s)-\Pi_{1,k}h(s)\right)'\exp\left(-\alpha\frac{t_k-s}{1-\alpha}\right)\mathrm{d}s\right| \nonumber\\
   &\leq\frac{\alpha \max\limits_{t\in[0,T]}h''(t)}{8(1-\alpha)^2}(k-1)\left(\Delta t\right)^3+\frac{\alpha \max\limits_{t\in[0,T]}h''(t)}{8(1-\alpha)^2}\left(\Delta t\right)^3\nonumber\\
   &=\frac{\alpha \max\limits_{t\in[0,T]}h''(t)}{8(1-\alpha)^2}k\left(\Delta t\right)^3.\nonumber
\end{align}
The estimate (\ref{lem-erro-Ft-eq1}) is proved. Hence
\begin{eqnarray}
   &&  \left|R_k\right|\leq\frac{\alpha\max\limits_{t\in[0,T]}h''(t)}{8(1-\alpha)^2}k\left(\Delta t\right)^3\leq\frac{\alpha T\max\limits_{t\in[0,T]}h''(t)}{8(1-\alpha)^2}\left(\Delta t\right)^2,\qquad j=1,2,\cdots M_t,\nonumber
\end{eqnarray}
which prove the conclusion of the lemma. \hfill$\Box$
\begin{rem}\label{rem-trun-err}
The second rate of convergence of L1 formula has been proven in \cite{akman2018new} by different methods, here we obtained identical results herein. Note that the rate of convergence of L1 formula for classical Caputo fractional derivative with order $\alpha$ is $2-\alpha$, this result seems reasonable since Caputo-Fabrizio derivative has smooth kernel.

\end{rem}

\subsection{Discretization in time}
We denote $u^k:=u(x,t_k)$ and $f^k(x):=f(x,t_k)$. Then from \eqref{app-CFD-rec-t1}-\eqref{app-CFD-rec-tk} and Lemma \ref{lem-erro-Ft}, the time fractional derivative \eqref{def-pCFD} at $t=t_k$ can be approximated by
\begin{eqnarray}
   && {}_0^{\mathrm{CF}}\!D_t^{\alpha_i}u^1(x)\approx F_t^{\alpha_i}u^1(x)=\frac{b_{1,1}^{(\alpha_i)}}{\alpha_i \Delta t}\left(u^1(x)-u^0(x)\right),\nonumber\\
   &&{}_0^{\mathrm{CF}}\!D_t^{\alpha_i}u^k(x)\approx F_t^{\alpha_i}u^k(x)=\exp\left(-\frac{\alpha_i\Delta t}{1-\alpha_i}\right)F_t^{\alpha_i}u^{k-1}(x)+\frac{b_{k,k}^{(\alpha_i)}}{\alpha_i \Delta t}\left(u^k(x)-u^{k-1}(x)\right),\quad k\geq 2,\nonumber
\end{eqnarray}
where
\begin{eqnarray}
   && b_{j,k}^{(\alpha_i)}:=\sigma_{j,k}^{(\alpha_i)}-\sigma_{j-1,k}^{(\alpha_i)},\quad  \sigma_{j,k}^{(\alpha_i)}:=\exp\left(-\alpha_i\frac{t_{k}-t_j}{1-\alpha_i}\right),\quad 1\leq k\leq N_T,\quad 1\leq j\leq k,\label{def-b}
\end{eqnarray}
 Then Eq.\eqref{pro} can be rewritten as
\begin{eqnarray}
   &&\sum_{i=1}^n\frac{d_i}{\alpha_i \Delta t}b_{1,1}^{(\alpha_i)}\left(u^1(x)-u^0(x)\right)=\frac{\partial^2u^1(x)}{\partial x^2}+f^1(x)+\sum_{i=1}^nd_iR_1^{(\alpha_i)},\nonumber\\
   && \sum_{i=1}^nd_i\exp\left(-\frac{\alpha_i\Delta t}{1-\alpha_i}\right)F_t^{\alpha_i}u^{k-1}(x)+\sum_{i=1}^n\frac{d_i}{\alpha_i \Delta t}b_{k,k}^{(\alpha_i)}\left(u^k(x)-u^{k-1}(x)\right)\nonumber\\
   &&\qquad\qquad=\frac{\partial^2u^k(x)}{\partial x^2}+f^k(x)+\sum_{i=1}^nd_iR_k^{(\alpha_i)},\qquad k\geq 2.\nonumber
\end{eqnarray}
where
\begin{eqnarray}
   && \left|R_k^{(\alpha_i)}\right|:=\left|{}_0^{\mathrm{CF}}\!D_t^{\alpha_i}u^k(x)-F_t^{\alpha_i} u^k(x)\right|\leq C_{u,\alpha_i}\left(\Delta t\right)^2,\label{err-cfd}
\end{eqnarray}
with $C_{u,\alpha_i}>0$ for $i=1,2,\cdots,n$. Notice that
\begin{eqnarray}
   && b_{k,k}^{(\alpha_i)}\equiv 1-\exp\left(-\frac{\alpha_i\Delta t}{1-\alpha_i}\right),\quad 1\leq k\leq N_T,\quad i=1,2,\cdots,n,\label{bkk}
\end{eqnarray}
we denote
\begin{eqnarray}
   && \kappa:=\eta_{k,k}^{-1},\qquad\eta_{l,s}:=\sum_{i=1}^n\frac{d_i}{\alpha_i \Delta t}b_{l,s}^{(\alpha_i)},\quad 1\leq s\leq N_T,\quad 1\leq l\leq k,\nonumber\\
   &&\widetilde{R}_{k}:=\kappa R_k,\qquad R_k:=\sum_{i=1}^nd_iR_k^{(\alpha_i)},\quad 1\leq k\leq N_T,\label{def-Rw}
\end{eqnarray}
the above equations are recast into
\begin{eqnarray}
   && u^1(x)-\kappa\frac{\partial^2u^1(x)}{\partial x^2}=u^{0}(x)+\kappa f^1(x)+\widetilde{R}_{1},\label{semi-R1}\label{Lt-sei-dis-R1}\\
   && u^k(x)-\kappa\frac{\partial^2u^k(x)}{\partial x^2}=u^{k-1}(x)-\kappa\sum_{i=1}^nd_i\exp\left(-\frac{\alpha_i\Delta t}{1-\alpha_i}\right)F_t^{\alpha_i}u^{k-1}(x)\nonumber\\
   &&\qquad\qquad\qquad\qquad\qquad+\kappa f^k(x)+\widetilde{R}_{k},\quad k\geq 2,\label{semi-Rk}\label{Lt-sei-dis-R2}
\end{eqnarray}
Let $u^k$ be the approximation for $u^k(x)$, and $f^k:=f^k(x)$. Then the semi-discrete problem of Eq. \eqref{pro} can be written as
\begin{eqnarray}
   &&  u^1-\kappa\frac{\partial^2u^1}{\partial x^2}=u^0+\kappa f^1, \label{Ft-sei-dis-1}\\
   &&  u^k-\kappa\frac{\partial^2u^k}{\partial x^2}=u^{k-1}-\kappa\sum_{i=1}^nd_i\exp\left(-\frac{\alpha_i\Delta t}{1-\alpha_i}\right)F_t^{\alpha_i}u^{k-1}+\kappa f^k,\quad k\geq 2,\label{Ft-sei-dis-3}\\
   &&u^0:=u(x,0)=\varphi(x),\quad x\in\Lambda,\label{Ft-sei-dis-5}\\
   &&u^k(0)=u^k(L)=0,\quad k=0,1,\cdots,N_T,\label{Ft-sei-dis-6}
\end{eqnarray}
where
\begin{eqnarray*}
   &&  F_t^{\alpha_i}u^1=\frac{1}{\alpha_i \Delta t}b_{1,1}^{(\alpha_i)}\left(u^1-u^0\right), \\
   && F_t^{\alpha_i}u^{k}=\exp\left(-\frac{\alpha_i\Delta t}{1-\alpha_i}\right)F_t^{\alpha_i}u^{k-1}+\frac{1}{\alpha_i \Delta t}b_{k,k}^{(\alpha_i)}\left(u^k-u^{k-1}\right),\quad k\geq 2.
\end{eqnarray*}
Moreover, by utilizing relation
\begin{eqnarray}
   && b_{j,k}^{(\alpha_i)}\exp\left(-\frac{\alpha_i\Delta t}{1-\alpha_i}\right)=b_{j,k+1}^{(\alpha_i)},\quad i=1,2,\cdots,n,\nonumber
\end{eqnarray}
we can easily derive an alternative formulation of \eqref{Ft-sei-dis-1}-\eqref{Ft-sei-dis-6} as follows
\begin{eqnarray}
   &&   u^1-\kappa\frac{\partial^2u^1}{\partial x^2}=u^0+\kappa f^1, \label{Lt-sei-dis-1}\\
   &&  u^k-\kappa\frac{\partial^2u^k}{\partial x^2}=\sum_{j=1}^{k-1}\left(\zeta_{j+1,k}-\zeta_{j,k}\right)u^j+\zeta_{1,k}u^0+\kappa f^k,\quad k\geq 2,\label{Lt-sei-dis-2}\\
   &&u^0:=u(x,0)=\varphi(x),\quad x\in\Lambda,\label{Lt-sei-dis-3}\\
   &&u^k(0)=u^k(L)=0,\quad k=0,1,\cdots,N_T,\label{Lt-sei-dis-4}
\end{eqnarray}
where
\begin{eqnarray}
   && \zeta_{j,k}:=\frac{\eta_{j,k}}{\eta_{k,k}},\quad 1\leq k\leq N_T,\quad 1\leq j\leq k. \label{def-zeta}
\end{eqnarray}
\begin{rem}
Since $L_t^{\alpha_i}=F_t^{\alpha_i}$, \eqref{Lt-sei-dis-1}-\eqref{Lt-sei-dis-4} can be also obtained by using $L_t^{\alpha_i}$ in Eq.\eqref{pro}.
It is noteworthy that equations \eqref{Ft-sei-dis-1}-\eqref{Ft-sei-dis-6} offer computational advantages over equations \eqref{Lt-sei-dis-1}-\eqref{Lt-sei-dis-4}. This is primarily attributed to the straightforward recurrence relation presented in equations \eqref{app-CFD-rec-t1} and \eqref{app-CFD-rec-tk}. However, \eqref{Lt-sei-dis-1}-\eqref{Lt-sei-dis-4} is more appropriate for our analysis than \eqref{Ft-sei-dis-1}-\eqref{Ft-sei-dis-6}, hence it play a crucial role in the subsequent sections.
\end{rem}
\begin{thm}\label{thm-Rk}
Let $\widetilde{R}_{k}$ be defined by \eqref{def-Rw}, then there exists a constant $\widetilde{c}>0$ such that
\begin{eqnarray}\label{err-Rk}
   &&  \left|\widetilde{R}_{k}\right|\leq \widetilde{c}\left(\Delta t\right)^{2},\quad k=1,2,\cdots,N_T.
\end{eqnarray}
\end{thm}
\noindent\emph{Proof.} Without loss of generality, we assume that $\Delta t\in(0,1)$. By the definition of $R_{k}$ and the inequalities of \eqref{err-cfd}, we have
\begin{eqnarray*}
 && \left|R_{k}\right|=\left|\sum_{i=1}^{n}d_iR_{k}^{\alpha_i}\right|\leq\sum_{i=1}^{n}d_i\left|R_{k}^{\alpha_i}\right|\leq\sum_{i=1}^{n}d_iC_{u,\alpha_i}(\Delta t)^{2}=\widehat{C}_{u,\alpha_i}\left(\Delta t\right)^{2},\quad k=1,2,\cdots,N_T.
\end{eqnarray*}
On the other hand, since
\begin{eqnarray*}
   &&\frac{1}{\alpha_i \Delta t} b_{k,k}^{(\alpha_i)}=\frac{1}{\alpha_i \Delta t}\left[1-\exp\left(-\frac{\alpha_i\Delta t}{1-\alpha_i}\right)\right]\to\frac{1}{1-\alpha_i},\qquad \Delta t\to 0,
\end{eqnarray*}
for $i=1,2,\cdots,n$, we get
\begin{eqnarray*}
   && \eta_{k,k}=\sum_{i=1}^n\frac{d_i}{\alpha_i \Delta t}b_{k,k}^{(\alpha_i)}\to\sum_{i=1}^n\frac{d_i}{1-\alpha_i},\qquad \Delta t\to 0.
\end{eqnarray*}
This implies that
\begin{eqnarray}
   && \left|\kappa\right|=\left|\eta_{k,k}^{-1}\right|=O(1).\label{kappa-bounded}
\end{eqnarray}
Therefore, there exists a constant $\widetilde{c}>0$ such that
\begin{eqnarray*}
  &&\left|\widetilde{R}_{k}\right|\leq|\kappa|\left|R_{k}\right|\leq \widetilde{c}\left(\Delta t\right)^{2},\quad k=1,2,\cdots,N_T,
\end{eqnarray*}
which prove the conclusion of the lemma. \hfill$\Box$

\begin{lem}\label{lem-bj}
Let the coefficients $b_{j,k}^{(\alpha_i)}$ be defined by \eqref{def-b}, then for every $i$,
\begin{eqnarray*}
   && 0<b_{j,N_T}^{(\alpha_i)}<\cdots<b_{j,k+1}^{(\alpha_i)}<b_{j,k}^{(\alpha_{i})}<\cdots<b_{j,j}^{(\alpha_{i})}<1, \\
   &&  0<b_{1,k}^{(\alpha_{i})}<\cdots<b_{j-1,k}^{(\alpha_{i})}< b_{j,k}^{(\alpha_{i})}<\cdots<b_{k,k}^{(\alpha_{i})}<1.
\end{eqnarray*}
\end{lem}
\noindent\emph{Proof.} $b_{j,k}^{(\alpha_i)}\in(0,1)$ can be easily obtained by the definition of $\sigma_{j,k}^{(\alpha_i)}$ and the monotone property of the function $g(x)=\exp(x)$. Finally, note that
\begin{eqnarray}
   &&  b_{j,k}^{(\alpha_i)}\exp\left(-\frac{\alpha_i\Delta t}{1-\alpha_i}\right)=b_{j,k+1}^{(\alpha_i)}=b_{j-1,k}^{(\alpha_i)},\quad i=1,2,\cdots,n.\label{rel-bjk-alpha}
\end{eqnarray}
Using the above equalities and the fact $\exp\left(-\frac{\alpha_i\Delta t}{1-\alpha_i}\right)\in(0,1)$ completes the proof of the lemma.\hfill$\Box$

\begin{rem}
\eqref{rel-bjk-alpha} gives a easy way to compute all the coefficients $b_{j,k}^{(\alpha_i)}$.
\end{rem}

\begin{lem}\label{lem-zetajk}
Let the coefficients $\zeta_{j,k}$ be defined by \eqref{def-zeta}, then
\begin{eqnarray*}
   && 0<\zeta_{j,N_T}<\cdots<\zeta_{j,k+1}<\zeta_{j,k}<\cdots<\zeta_{j,j}=1 \\
   &&  0<\zeta_{1,k}<\cdots<\zeta_{j-1,k}< \zeta_{j,k}<\cdots<\zeta_{k,k}=1.
\end{eqnarray*}
\end{lem}

\noindent\emph{Proof.} By Lemma \ref{lem-bj}, and the definition of $\zeta_{j,k}$, we can readily arrive at these conclusions.\hfill$\Box$

\subsection{Stability and convergence analysis of the semi-discrete scheme}
To discuss the stability and convergence of the semi-discrete scheme, we introduce functional spaces equipped with standard norms and inner products that will be utilized subsequently. Let $L^2(\Lambda)$ is the space of measurable functions whose square is Lebesgue integrable in $\Lambda$. Then
\begin{eqnarray*}
   && H^1(\Lambda):=\left\{v\in L^2(\Lambda),~\frac{\mathrm{d} v}{\mathrm{d} x}\in L^2(\Lambda)\right\},\\
   && H_0^1(\Lambda):=\left\{v\in H^1(\Omega),~v|_{\partial\Lambda}=0\right\},\\
   && H^m(\Lambda):=\left\{v\in L^2(\Lambda),~\frac{\mathrm{d}^k v}{\mathrm{d} x^k}\in L^2(\Lambda),~\mathrm{for~all~positive~integer}~k\leq m\right\},
\end{eqnarray*}
The inner products of $L^2(\Lambda)$ and $H^1(\Lambda)$ are defined, respectively, by
\begin{eqnarray*}
   && (u,v)=\int_{\Lambda}uv~\mathrm{d}x,\qquad (u,v)_1=(u,v)+\left(\frac{\mathrm{d} u}{\mathrm{d} x},\frac{\mathrm{d} v}{\mathrm{d} x}\right),
\end{eqnarray*}
and the corresponding norms by
\begin{eqnarray*}
   && \|v\|_0=(v,v)^{1/2},\qquad \|v\|_1=(v,v)_1^{1/2}.
\end{eqnarray*}
The norm $\|\cdot\|_m$ of the space $H^m(\Omega)$~$(m\in\mathrm{N})$ is defined by
\begin{eqnarray*}
   && \|v\|_m:=\left(\sum\limits_{k=0}^m\left\|\frac{\mathrm{d}^k v}{\mathrm{d} x^k}\right\|_0^2\right)^{\frac{1}{2}}.
\end{eqnarray*}
In this paper, instead of using the above standard $H^1$-norm, we prefer to define $\|\cdot\|_1$ by
\begin{eqnarray}
   && \|v\|_1=\left(\|v\|_0^2+\kappa\left\|\frac{\mathrm{d} v}{\mathrm{d} x}\right\|_0^2\right)^{1/2}.\label{new-H1-norm}
\end{eqnarray}
It is widely acknowledged that the standard $H^1$-norm and the norm defined by \eqref{new-H1-norm} are equivalent; therefore, we will adopt the latter in subsequent discussions.

The variational (weak) formulation of the Eqs.\eqref{Ft-sei-dis-1} and \eqref{Ft-sei-dis-3}/\eqref{Lt-sei-dis-2}, subject to the boundary condition \eqref{Ft-sei-dis-6}, can be expressed as finding $u^{k}\in H_0^1(\Lambda)$ such that for $\forall v\in H_0^1(\Lambda)$
\begin{align}
\left(u^{1},v\right)+\kappa\left(\frac{\partial u^1}{\partial x},\frac{\partial v}{\partial x}\right)&=\left(u^{0},v\right)+\kappa\left(f^{1},v\right),\label{var-sem-dis-k0}\\
    \left(u^{k},v\right)+\kappa\left(\frac{\partial u^{k}}{\partial x},\frac{\partial v}{\partial x}\right)&=
   \left(u^{k-1},v\right)-\eta_{k,k}^{-1}\sum_{i=1}^nd_i\exp\left(-\frac{\alpha_i\Delta t}{1-\alpha_i}\right)\left(F_t^{\alpha_i}u^{k-1},v\right)+\kappa \left(f^k,v\right),
    \nonumber\\
 &=\sum_{j=1}^{k-1}\left(\zeta_{j+1,k}-\zeta_{j,k}\right)\left(u^j,v\right)+\zeta_{1,k}\left(u^0,v\right)+\kappa \left(f^k,v\right), \quad k\geq2, \label{var-sem-dis-k1}
\end{align}
where $u^0:=u(x,0)$.

For the semi-discretized problem \eqref{var-sem-dis-k0}-\eqref{var-sem-dis-k1}, we can establish a stability result as follows.
\begin{thm}\label{thm-stable-sem-dis}
The semi-discretized problem \eqref{var-sem-dis-k0}-\eqref{var-sem-dis-k1} is unconditionally stable in the sense that for all $\Delta t>0$, it holds
\begin{eqnarray*}
   && \|u^{k}\|_1\leq\|u^0\|_0+\frac{\kappa}{\zeta_{1,k}}\max_{1\leq l\leq N_T}\|f^l\|_0,\quad k=1,2,\cdots,N_T.
\end{eqnarray*}
\end{thm}
\noindent\emph{Proof.} By mathematical induction. First of all, when $k=1$, we have
\begin{eqnarray*}
   && \left(u^1,v\right)+\kappa\left(\frac{\partial u^1}{\partial x},\frac{\partial v}{\partial x}\right)=\left(u^0,v\right)+\kappa\left(f^{1},v\right),\quad \forall v\in H_0^1(\Lambda).
\end{eqnarray*}
Notice that $\|v\|_0\leq\|v\|_1$, taking $v=u^1$ and using the Cauchy-Schwarz inequality, we obtain immediately
\begin{eqnarray*}
   && \|u^1\|_1\leq\|u^0\|_0+\kappa\|f^1\|_0=\|u^0\|_0+\frac{\kappa}{\zeta_{1,1}}\|f^1\|_0\leq\|u^0\|_0+\frac{\kappa}{\zeta_{1,1}}\max_{1\leq l\leq N_T}\|f^l\|_0.
\end{eqnarray*}
Now, suppose
\begin{eqnarray}
   && \|u^s\|_1\leq\|u^0\|_0+\frac{\kappa}{\zeta_{1,s}}\max_{1\leq l\leq N_T}\|f^l\|_0,\quad s=1,2,\cdots,k-1.\label{assum1}
\end{eqnarray}
Taking $v=u^{k}$ in (\ref{var-sem-dis-k1}) gives
\begin{eqnarray*}
   && \|u^{k}\|_1^2\leq  \left[\sum_{j=1}^{k-1}\left(\zeta_{j+1,k}-\zeta_{j,k}\right)\|u^{j}\|_0+\zeta_{1,k}\|u^0\|_0+\kappa \|f^{k}\|_0\right]\|u^{k}\|_0.
\end{eqnarray*}
Hence, by using (\ref{assum1}) and Lemma \ref{lem-zetajk}, we have
\begin{align*}
   \|u^{k}\|_1 & \leq\sum_{j=1}^{k-1}\left(\zeta_{j+1,k}-\zeta_{j,k}\right)\|u^{j}\|_0+\zeta_{1,k}\|u^0\|_0+\kappa \|f^{k}\|_0\\
   &\leq \left[\sum_{j=1}^{k-1}\left(\zeta_{j+1,k}-\zeta_{j,k}\right)
   +\zeta_{1,k}\right]\|u^0\|_0+ \left[\sum_{j=1}^{k-1}\frac{\zeta_{j+1,k}-\zeta_{j,k}}{\zeta_{1,j}}+1\right]\kappa\max_{1\leq l\leq N_T}\|f^l\|_0\\
   &\leq\|u^0\|_0+\left[\sum_{j=1}^{k-1}\frac{\zeta_{j+1,k}-\zeta_{j,k}}{\zeta_{1,k}}+1\right]\kappa\max_{1\leq l\leq N_T}\|f^l\|_0\\
   &=\|u^0\|_0+\frac{\kappa}{\zeta_{1,k}}\max_{1\leq l\leq N_T}\|f^l\|_0.
\end{align*}
Thus, the proof is completed.\hfill$\Box$
\begin{rem}
In the proof of the following theorem, we will demonstrate that $\zeta_{1,k}^{-1}$ is bounded. As shown in equation \eqref{kappa-bounded}, $|\kappa|=O(1)$. Therefore, it follows that $\kappa\zeta_{1,k}^{-1}$  is also bounded.
\end{rem}

 We now conduct an error analysis for the solution of the semi-discretized problem.
 \begin{thm}\label{thm-sem-err}
 Assuming $n\neq 1$ and $0<\alpha_1\leq\alpha_2\leq\cdots\leq\alpha_{n-1}\leq\alpha_n<1$.
Let $u^k(x)$ be the exact solution of \eqref{pro}-\eqref{pro-bou-cond}, $\{u^k\}_{k=0}^{N_T}$ be the solution of semi-discretized problem \eqref{var-sem-dis-k0}-\eqref{var-sem-dis-k1} with initial condition $u^0:=u^0(x)$, then the following error estimates hold:
\begin{eqnarray}
   && \|u^k(x)-u^k\|_1\leq \widetilde{c}S\exp\left(\frac{\alpha_nT}{1-\alpha_n}\right)\left(\Delta t\right)^{2},\quad k=1,2,\cdots,N_T,\label{thm-sem-err-1}
\end{eqnarray}
where the constant $\widetilde{c}$ is defined in \eqref{err-Rk} and $S$ is the length of $\Lambda$.
\end{thm}

\noindent\emph{Proof.} We shall establish the following estimate through a process of mathematical induction:
\begin{equation}\label{thm-sem-err-m1}
    \|u^k(x)-u^k\|_1\leq \frac{\widetilde{c}S}{\zeta_{1,k}}\left(\Delta t\right)^{2},\quad k=1,2,\cdots,N_T.
\end{equation}
Let $\bar{e}^k=u^k(x)-u^k$,~$k=1,2,\cdots,N_t$. By combining \eqref{Lt-sei-dis-R1} and \eqref{var-sem-dis-k0}, the error $\bar{e}^1$ satisfies
\begin{eqnarray*}
   &&  \left(\bar{e}^1,v\right)+\kappa\left(\frac{\partial \bar{e}^1}{\partial x},\frac{\partial v}{\partial x}\right)=\left(\bar{e}^0,v\right)+\left(\widetilde{R}_{1},v\right)=\left(\widetilde{R}_{1},v\right),\quad \forall v\in H_0^1(\Lambda).
\end{eqnarray*}
Taking $v=\bar{e}^1$ yields $\|\bar{e}^1\|_1^2\leq\|\widetilde{R}_{1}\|_0\|\bar{e}^1\|_0$. This, in conjunction with \eqref{err-Rk}, yields
\begin{eqnarray*}
   && \|u^1(x)-u^1\|_1=\|\bar{e}^1\|_1\leq\|\widetilde{R}_{1}\|_0\leq \widetilde{c}S\left(\Delta t\right)^{2}=\frac{\widetilde{c}S}{\zeta_{1,1}}\left(\Delta t\right)^{2}.
\end{eqnarray*}
Therefore, \eqref{thm-sem-err-m1} holds for $k=1$. Assuming that \eqref{thm-sem-err-m1} holds for all $k=1,2,\cdots,l-1$, it is necessary to demonstrate its validity for $k=l$. By combining \eqref{Lt-sei-dis-R1}, \eqref{Lt-sei-dis-R2} and \eqref{var-sem-dis-k1}, for $\forall v\in H_0^1(\Lambda)$, we have
\begin{eqnarray*}
   &&   \left(\bar{e}^{l},v\right)+\kappa\left(\frac{\partial \bar{e}^{l}}{\partial x},\frac{\partial v}{\partial x}\right)=\sum_{j=1}^{l-1}\left(\zeta_{j+1,l}-\zeta_{j,l}\right)\left(\bar{e}^j,v\right)+\zeta_{1,l}\left(\bar{e}^0,v\right)+\left(\widetilde{R}_{l },v\right).\nonumber
\end{eqnarray*}
Let $v=\bar{e}^{l}$ in the above equation, then
\begin{eqnarray*}
   && \|\bar{e}^{l}\|_1\leq \sum_{j=1}^{l-1}\left(\zeta_{j+1,l}-\zeta_{j,l}\right)\|\bar{e}^{j}\|_0+\zeta_{1,l}\|\bar{e}^{0}\|_0+\|\widetilde{R}_{l}\|_0.
\end{eqnarray*}
Using the induction assumption and Lemma \ref{lem-zetajk}, we derive
\begin{eqnarray*}
   && \|\bar{e}^{l}\|_1\leq\sum_{j=1}^{l-1}\frac{\zeta_{j+1,l}-\zeta_{j,l}}{\zeta_{1,j}}\widetilde{c}S\left(\Delta t\right)^{2}+\widetilde{c}S\left(\Delta t\right)^{2}\leq\left[\sum_{j=1}^{l-1}\frac{\zeta_{j+1,l}-\zeta_{j,l}}{\zeta_{1,l}}+1\right]\widetilde{c}S\left(\Delta t\right)^{2}=\frac{\widetilde{c}S}{\zeta_{1,l}}\left(\Delta t\right)^{2}.
\end{eqnarray*}

Next we show that $\zeta_{1,k}^{-1}$ is bounded. Considering that function $h(\alpha)=:-\frac{\alpha}{1-\alpha}$ is decreasing on $(0,1)$, we have
\begin{eqnarray*}
   && b_{1,k}^{(\alpha_i)}=\exp\left(-\frac{\alpha_i(k-1)\Delta t}{1-\alpha_i}\right)b_{k,k}^{(\alpha_i)}\geq \exp\left(-\frac{\alpha_n(k-1)\Delta t}{1-\alpha_n}\right)b_{k,k}^{(\alpha_i)},\quad k=1,2,\cdots,N_T.
\end{eqnarray*}
by combining equation \eqref{rel-bjk-alpha}.
Therefore,
\begin{eqnarray*}
   && \frac{1}{\zeta_{1,k}}=\frac{\eta_{k,k}}{\eta_{1,k}}=\frac{\sum\limits_{i=1}^n\frac{d_i}{\alpha_i \Delta t}b_{k,k}^{(\alpha_i)}}{\sum\limits_{i=1}^n\frac{d_i}{\alpha_i \Delta t}b_{1,k}^{(\alpha_i)}}\leq \exp\left(\frac{\alpha_n(k-1)\Delta t}{1-\alpha_n}\right)\leq\exp\left(\frac{\alpha_nT}{1-\alpha_n}\right),\quad k=1,2,\cdots,N_T.
\end{eqnarray*}
Consequently we obtain, for all $k$ such that $k\Delta t\leq T$,
\begin{eqnarray*}
   &&  \|u^k(x)-u^k\|_1\leq \widetilde{c}S\exp\left(\frac{\alpha_nT}{1-\alpha_n}\right)\left(\Delta t\right)^{2},\quad k=1,2,\cdots,N_T.
\end{eqnarray*}
\hfill$\Box$

\section{Full discretization}
\subsection{A shifted Legendre collocation method in space}
We shall begin by providing a comprehensive overview of fundamental definitions and properties pertaining to Legendre Gauss-type quadratures.
Let $\mathds{P}_N(\Lambda)$ denote the space of algebraic polynomials of degree less than or
equal to $N$ with respect to variable $x$, and $L_N(x)$ be the Legendre polynomial
of degree $N$ on the interval $[-1,1]$. Then the discrete space, denoted by $\mathds{P}_N^0(\Lambda):=\mathds{P}_N(\Lambda)\cap H_0^1(\Lambda)$.

Let $\pi_N^1$ be the $H^1$-orthogonal projection operator from $H_0^1(\Lambda)$ into $\mathds{P}_N^0(\Lambda)$, associated to the norm $\|\cdot\|_1$ defined in (\ref{new-H1-norm}), that is, for all $\psi\in H_0^1(\Lambda)$, define $\pi_N^1\psi\in\mathds{P}_N^0(\Lambda)$, such that, $\forall v_N\in\mathds{P}_N^0(\Lambda)$,
\begin{equation}\label{h1-orth-proj}
    \left(\pi_N^1\psi,v_N\right)+\kappa\left(\frac{\mathrm{d}}{\mathrm{d} x}\pi_N^1\psi,\frac{\mathrm{d}}{\mathrm{d} x} v_N\right)=\left(\psi,v_N\right)+\kappa\left(\frac{\mathrm{d}}{\mathrm{d} x}\psi,\frac{\mathrm{d}}{\mathrm{d} x} v_N\right).
\end{equation}
From \cite{bernardi1992approximations}, the following estimate of projection holds:
\begin{equation}\label{pro-est}
   \|\psi-\pi_N^1\psi\|_1\leq c_1N^{1-m}\|\psi\|_m,\quad \forall \psi\in H^m(\Lambda)\cap H_0^1(\Lambda),\quad m\geq 1.
\end{equation}

 Define the Legendre-Gauss-Lobatto nodes and weights as $\xi_p$ and $\omega_p$, $p=0,1,\cdots,N$, $N\geq 1$,
where $\{\xi_k\}_{k=0}^{N}$ are the zeroes of $(1-x^2)L'_N(x)$, and
\begin{eqnarray*}
   && \omega_k=\frac{2}{N(N+1)}\frac{1}{L_N^2(\xi_k)},\quad k=0,1,\cdots,N.
\end{eqnarray*}
Moreover, the following quadrature holds
\begin{eqnarray*}
   && \int_{-1}^1\phi(x)\mathrm{d}x=\sum_{k=0}^N\phi(\xi_k)\omega_k,\quad \forall \phi(x)\in \mathds{P}_{2N-1}([-1,1]).
\end{eqnarray*}
The discrete inner product and norm defined as follow, for any continuous functions $\phi,\psi\in \mathds{C}\left([-1,1]\right)$,
\begin{eqnarray*}
   && \left(\phi,\psi\right)_N:=\sum\limits_{p=0}^{N}\phi(\xi_{p})\psi(\xi_{p})\omega_{p},\qquad \|\phi\|_N:=(\phi,\phi)_N^{1/2},
\end{eqnarray*}
From \cite{quarteroni2008numerical}, the discrete norm $\|\cdot\|_N$ is equivalent to the standard $L^2$-norm in $\mathds{P}_N([-1,1])$. If we denote $\{\hat{\xi}_{p}\}_{p=0}^N$
and $\{\hat{\omega}_{p}\}_{p=0}^N$ as the nodes and
weights of shifted Legendre-Gauss-Lobatto quadratures on $\overline{\Lambda}$, then one can easily show that
\begin{eqnarray*}
   && \hat{\xi}_k=\frac{S}{2}\left(\xi_k+1\right),\quad \hat{\omega}_k=\frac{S}{2}\omega_k,\quad k=0,1,\cdots,N.
\end{eqnarray*}
Thus, we define the discrete inner product and norm on $\overline{\Lambda}$ as follows
\begin{eqnarray*}
   &&  \left(\phi,\psi\right)_{\widehat{N}}:=\sum\limits_{p=0}^{N}\phi(\hat{\xi}_{p})\psi(\hat{\xi}_{p})\hat{\omega}_{p},\quad \|\phi\|_{\widehat{N}}:=(\phi,\phi)_{\widehat{N}}^{1/2},\quad \forall\phi,\psi\in \mathds{C}(\overline{\Lambda}).
\end{eqnarray*}
It is not  difficult to obtain that
\begin{eqnarray}
   && \|u_N\|_0\leq\|u_N\|_{\widehat{N}}\leq\sqrt{3}\|u_N\|_0,\quad \forall u_N\in\mathds{P}_N(\Lambda), \label{eqi-norm}
\end{eqnarray}
and
\begin{equation}\label{app-norm-normN}
     \left(\phi,\psi\right)_{\widehat{N}}= \left(\phi,\psi\right),\quad \forall \phi\psi\in\mathds{P}_{2N-1}(\Lambda).
\end{equation}
We introduce the operator of interpolation at the $N+1$ shifted Legendre-Gauss-Lobatto nodes, denoted by $I_N$, i.e., $\forall \psi\in C(\overline{\Lambda})$,~$I_N\psi\in\mathds{P}_N(\Lambda)$, such that
\begin{eqnarray}
   && I_N\psi(\hat{\xi}_p)=\psi(\hat{\xi}_p),\quad p=0,1,\cdots,N.\label{inter-cond}
\end{eqnarray}
The interpolation error estimate (see \cite{bernardi1992approximations}) is
\begin{equation}\label{err-inter}
  \|\psi-I_N\psi\|_1\leq c_2N^{1-m}\|\psi\|_m,\quad\forall\psi\in H^m(\Lambda)\cap H_0^1(\Lambda),\quad m\geq 1.
\end{equation}

Now consider the spectral discretization to the problem $(\ref{var-sem-dis-k0})$-$(\ref{var-sem-dis-k1})$ as follows: find $u_N^{k}\in\mathds{P}_N^0(\Lambda)$, such that
\begin{eqnarray}
    &&A_{\widehat{N}}(u_N^{k},v_N)=F_{\widehat{N}}^{k}(v_N), \quad\forall v_N\in\mathds{P}_N^0(\Lambda),\quad k\geq 1,\label{var-full-dis-k1-coll}
\end{eqnarray}
where
\begin{eqnarray*}
&& A_{\widehat{N}}(u_N^{k},v_N):=\left(u_N^{k},v_N\right)_{\widehat{N}}+\kappa\left(\frac{\partial u_N^{k}}{\partial x},\frac{\partial v_N}{\partial x}\right)_{\widehat{N}},\quad k\geq 1,\nonumber\\
  && F_{\widehat{N}}^1(v_N):=\left(u_N^{0},v_N\right)_{\widehat{N}}+\kappa\left(f^{1},v_N\right)_{\widehat{N}},\nonumber\\
   &&F_{\widehat{N}}^{k}(v_N):=\left(u_N^{k-1},v_N\right)_{\widehat{N}}-\eta_{k,k}^{-1}\sum_{i=1}^nd_i\exp\left(-\frac{\alpha_i\Delta t}{1-\alpha_i}\right)\left(F_t^{\alpha_i}u_N^{k-1},v_N\right)_{\widehat{N}}+\kappa \left(f^k,v_N\right)_{\widehat{N}},
    \nonumber\\
  &&\qquad\qquad= \sum_{j=1}^{k-1}\left(\zeta_{j+1,k}-\zeta_{j,k}\right)\left(u_N^j,v_N\right)_{\widehat{N}}+\zeta_{1,k}\left(u_N^0,v_N\right)_{\widehat{N}}+\kappa \left(f^k,v_N\right)_{\widehat{N}},\quad k\geq 2,\nonumber\\
   &&u_N^0:=I_Nu^0(x).\nonumber
\end{eqnarray*}
For $\{u_N^j\}_{j=0}^{k-1}$ given, the well-posedness of the problem \eqref{var-full-dis-k1-coll} is guaranteed by the well-known Lax-Milgram Lemma.

\subsection{Convergence analysis of the full discretization scheme}
To simplify matters, we present the semi-discretized problem \eqref{var-sem-dis-k0}-\eqref{var-sem-dis-k1} in a compact form: find $u^{k}\in H_0^1(\Lambda)$, such that
\begin{eqnarray*}
   && A\left(u^{k},v\right)=F^{k}\left(v\right),\quad  \forall v\in H_0^1(\Lambda),\quad k\geq 1,
\end{eqnarray*}
where
\begin{eqnarray*}
&& A(u^{k},v):=\left(u^{k},v\right)+\kappa\left(\frac{\partial u^{k}}{\partial x},\frac{\partial v}{\partial x}\right),\quad k\geq 1,\nonumber\\
  && F^1(v):=\left(u^{0},v\right)+\kappa\left(f^{1},v\right),\nonumber\\
   &&F^{k}(v):=\left(u^{k-1},v\right)-\eta_{k,k}^{-1}\sum_{i=1}^nd_i\exp\left(-\frac{\alpha_i\Delta t}{1-\alpha_i}\right)\left(F_t^{\alpha_i}u^{k-1},v\right)+\kappa \left(f^k,v\right),
    \nonumber\\
  &&\qquad\quad= \sum_{j=1}^{k-1}\left(\zeta_{j+1,k}-\zeta_{j,k}\right)\left(u^j,v\right)+\zeta_{1,k}\left(u^0,v\right)+\kappa \left(f^k,v\right),\quad k\geq 2,\nonumber\\
   &&u^0:=u^0(x).\nonumber
\end{eqnarray*}
We denote by $\|\cdot\|_{1,\widehat{N}}$ the norm associated to the bilinear form $A_{\widehat{N}}(\cdot,\cdot)$:
\begin{eqnarray*}
   && \|\phi\|_{1,\widehat{N}}:=A_{\widehat{N}}^{1/2}(\phi,\phi),\quad \forall \phi\in \mathds{C}(\overline{\Lambda}).
\end{eqnarray*}
It follows from \eqref{eqi-norm} that for all $v_N\in\mathds{P}_N(\Lambda)$ the
discrete norm $\|\cdot\|_{1,\widehat{N}}$ is equivalent to the norm $\|\cdot\|_1$ defined in \eqref{new-H1-norm}.

\begin{thm}\label{thm-full}
Assuming $n\neq 1$ and $0<\alpha_1\leq\alpha_2\leq\cdots\leq\alpha_{n-1}\leq\alpha_n<1$.
Let $\{u_N^k\}_{k=0}^{N_T}$ is the solution of the problem \eqref{var-full-dis-k1-coll} with the initial condition $u_N^0$ taken to be $I_Nu^0(x)$, $\{u^k\}_{k=0}^{N_T}$ the solution of the semi-discretized problem \eqref{var-sem-dis-k0}-\eqref{var-sem-dis-k1}. Suppose that $u^k\in H^m(\Lambda)\cap H_0^1(\Lambda)$ with $m> 1$, for $k=1,2,\cdots,N_T$,
then there exists a constant $\overline{c}>0$ such that
\begin{eqnarray*}
   && \|u^k-u_N^{k}\|_{1,{\widehat{N}}} \leq
  \overline{c}\exp\left(\frac{\alpha_nT}{1-\alpha_n}\right)\left( N^{-m}\max\limits_{0\leq l\leq k}\|f^l\|_m+(N-1)^{1-m}\max\limits_{0\leq l\leq k}\|u^l\|_m\right);
\end{eqnarray*}

\end{thm}
\noindent\emph{Proof.} For any $\forall v_{N-1}\in\mathds{P}_{N-1}^0(\Lambda)$, denote $\rho_N^{k}:=u_N^{k}-v_{N-1}$. It is direct to check that
\begin{eqnarray*}
   && A_{\widehat{N}}\left(\rho_N^{k},\rho_N^{k}\right)=A\left(u^{k}-v_{N-1},\rho_N^{k}\right)+A\left(v_{N-1},\rho_N^{k}\right)
   -A_{\widehat{N}}\left(v_{N-1},\rho_N^{k}\right)+F_{\widehat{N}}^{k}\left(\rho_N^{k}\right)-F^{k}\left(\rho_N^{k}\right).
\end{eqnarray*}
By virtue of (\ref{app-norm-normN}) gives
\begin{eqnarray*}
   && A\left(v_{N-1},\rho_N^{k}\right)=A_{\widehat{N}}\left(v_{N-1},\rho_N^{k}\right),\quad \forall v_{N-1}\in\mathrm{P}_{N-1}^0(\Lambda),
\end{eqnarray*}
hence
\begin{equation}\label{thm-col-eq1}
   \|\rho_N^{k}\|_{1,{\widehat{N}}}^2\leq\|u^{k}-v_{N-1}\|_1\|\rho_N^{k}\|_1+|F^{k}(\rho_N^{k})-F_{\widehat{N}}^{k}(\rho_{N}^{k})|,\quad \forall v_{N-1}\in\mathds{P}_{N-1}^0(\Lambda).
\end{equation}
For the last term, by definition, we have
\begin{equation}\label{thm-col-eq2}
  F^{1}(\rho_N^{1})-F_{\widehat{N}}^{1}(\rho_N^{1})=\left[\left(u^0,\rho_N^{1}\right)-\left(u_N^0,\rho_N^{1}\right)_{\widehat{N}}\right]+\kappa\left[\left(f^1,\rho_N^{1}\right)-\left(f^1,\rho_N^{1}\right)_{\widehat{N}}\right],
\end{equation}
and
\begin{eqnarray}
   &&  F^{k}(\rho_N^{k})-F_{\widehat{N}}^{k}(\rho_N^{k})=\sum_{j=1}^{k-1}\left(\zeta_{j+1,k}-\zeta_{j,k}\right)\left[\left(u^j,\rho_N^k\right)-\left(u_N^j,\rho_N^{k}\right)_{\widehat{N}}\right]
   \nonumber  \\ &&\qquad\qquad+\zeta_{1,k}\left[\left(u^0,\rho_N^k\right)-\left(u_N^0,\rho_N^{k}\right)_{\widehat{N}}\right]+\left[\left(f^{k},\rho_N^{k}\right)-\left(f^{k},\rho_N^{k}\right)_{\widehat{N}}\right],\quad k\geq 2.  \label{thm-col-eq3}
\end{eqnarray}
It is known that the following result holds (see e.g. \cite{quarteroni2008numerical,bernardi1992approximations}): $\forall g\in H^m(\Lambda)$, $m\geq 1$, $\forall \delta_N\in \mathds{P}_N(\Lambda)$,
\begin{eqnarray*}
   && \left|\left(g,\delta_N\right)-\left(g,\delta_N\right)_{\widehat{N}}\right|\leq c_3N^{-m}\|g\|_m\|\delta_N\|_0.
\end{eqnarray*}
Thus for $\forall g\in H^m(\Lambda)$, $m\geq 1$, $\forall g_N,\rho_N\in\mathds{P}_N(\Lambda)$ we have
\begin{align*}
  \left| \left(g,\rho_N\right)-\left(g_N,\rho_N\right)_{\widehat{N}}\right|&\leq \left| \left(g,\rho_N\right)-\left(g,\rho_N\right)_{\widehat{N}}\right|
   +\left| \left(g,\rho_N\right)_{\widehat{N}}-\left(g_N,\rho_N\right)_{\widehat{N}}\right| \\
   &\leq c_3N^{-m}\|g\|_m\|\rho_N\|_0+\|g-g_N\|_{\widehat{N}}\|\rho_N\|_{\widehat{N}}\\
   &\leq \left(c_3N^{-m}\|g\|_m+\|g-g_N\|_{\widehat{N}}\right)\|\rho_N\|_{\widehat{N}}.
\end{align*}
Applying the above results to (\ref{thm-col-eq2}) and (\ref{thm-col-eq3}), we obtain
\begin{eqnarray*}
   && \left|F^{1}(\rho_N^{1})-F_{\widehat{N}}^{1}(\rho_N^{1})\right|\leq \left(c_3N^{-m}\|u^0\|_m+\|u^0-u_N^0\|_{\widehat{N}}+c_3\kappa N^{-m}\|f^1\|_m\right)\|\rho_N^{1}\|_{\widehat{N}},
\end{eqnarray*}
and
\begin{eqnarray*}
   && \left|F^{k}(\rho_N^{k})-F_{\widehat{N}}^{k}(\rho_N^{k})\right|\leq\bigg[\sum\limits_{j=1}^{k-1}\left(\zeta_{j+1,k}-\zeta_{j,k}\right)\|u^j-u_N^{j}\|_{\widehat{N}}+\zeta_{1,k}\|u^0-u_N^0\|_{\widehat{N}}\\
&& \qquad\qquad+c_3N^{-m}\max\limits_{0\leq j\leq k}\|u^j\|_m+c_3\kappa N^{-m}\|f^{k}\|_m\bigg]\|\rho_N^{k}\|_{\widehat{N}},\quad k\geq 2.
\end{eqnarray*}
Let $\varepsilon_N^j:=u^j-u_N^j$, using (\ref{thm-col-eq1}) and the norm equivalence, for $\forall v_{N-1}\in\mathds{P}_{N-1}^0(\Lambda)$, we have
\begin{eqnarray*}
   &&   \|\rho_N^{1}\|_{1,{\widehat{N}}}\leq\left\|\varepsilon_N^0\right\|_{\widehat{N}}+c_3N^{-m}\left\|u^0\right\|_m+c_3\kappa N^{-m}\left\|f^1\right\|_m+c_4\|u^1-v_{N-1}\|_{1,{\widehat{N}}},
\end{eqnarray*}
and
\begin{eqnarray*}
   &&  \|\rho_N^{k}\|_{1,{\widehat{N}}}\leq \sum\limits_{j=1}^{k-1}\left(\zeta_{j+1,k}-\zeta_{j,k}\right)\|\varepsilon_N^j\|_{\widehat{N}}+\zeta_{1,k}\|\varepsilon_N^0\|_{\widehat{N}}
   +c_3N^{-m}\max\limits_{0\leq j\leq k}\|u^j\|_m\\
   &&\qquad\qquad\quad+c_3\kappa N^{-m}\|f^{k}\|_m+c_4\|u^k-v_{N-1}\|_{1,{\widehat{N}}},\quad k\geq 2.
\end{eqnarray*}
By triangular inequality
\begin{eqnarray*}
   && \|\varepsilon_N^{k}\|_{1,{\widehat{N}}}\leq\|\rho_N^{k}\|_{1,{\widehat{N}}}+\|u^{k}-v_{N-1}\|_{1,{\widehat{N}}},
\end{eqnarray*}
for $\forall v_{N-1}\in\mathds{P}_{N-1}^0(\Lambda)$, we obtain
\begin{eqnarray*}
   &&   \|\varepsilon_N^{1}\|_{1,{\widehat{N}}}\leq\|\varepsilon_N^0\|_{\widehat{N}}+c_3N^{-m}\|u^0\|_m+c_3\kappa N^{-m}\|f^1\|_m+c_5\|u^1-v_{N-1}\|_{1,{\widehat{N}}},
\end{eqnarray*}
and
\begin{eqnarray*}
   &&  \|\varepsilon_N^{k}\|_{1,{\widehat{N}}}\leq \sum\limits_{j=1}^{k-1}\left(\zeta_{j+1,k}-\zeta_{j,k}\right)\|\varepsilon_N^j\|_{\widehat{N}}+\zeta_{1,k}\|\varepsilon_N^0\|_{\widehat{N}}
   +c_3N^{-m}\max\limits_{0\leq j\leq k}\|u^j\|_m\\
   &&\qquad\qquad\quad+c_3\kappa N^{-m}\|f^{k}\|_m+c_5\|u^k-v_{N-1}\|_{1,{\widehat{N}}},\quad k\geq 2.
\end{eqnarray*}
The above estimate specially holds for $v_{N-1}=\pi_{N-1}^1u^{k}$, which implies
\begin{align*}
      \|\varepsilon_N^{1}\|_{1,{\widehat{N}}}&\leq\|\varepsilon_N^0\|_{\widehat{N}}+c_3N^{-m}\|u^0\|_m+c_3\kappa N^{-m}\|f^1\|_m+c_6(N-1)^{1-m}\|u^1\|_m\\
    &\leq\|\varepsilon_N^0\|_{\widehat{N}}+c_3\kappa N^{-m}\max\limits_{0\leq l\leq 1}\|f^l\|_m+c_7(N-1)^{1-m}\max\limits_{0\leq l\leq 1}\|u^l\|_m,
\end{align*}
and
\begin{align*}
    \|\varepsilon_N^{k}\|_{1,{\widehat{N}}}&\leq \sum\limits_{j=1}^{k-1}\left(\zeta_{j+1,k}-\zeta_{j,k}\right)\|\varepsilon_N^j\|_{\widehat{N}}+\zeta_{1,k}\|\varepsilon_N^0\|_{\widehat{N}}
   +c_3N^{-m}\max\limits_{0\leq j\leq k}\|u^j\|_m\\
   &\quad +c_3\kappa N^{-m}\|f^{k}\|_m+c_6(N-1)^{1-m}\|u^k\|_m\\
   &\leq \sum\limits_{j=1}^{k-1}\left(\zeta_{j+1,k}-\zeta_{j,k}\right)\|\varepsilon_N^j\|_{\widehat{N}}+\zeta_{1,k}\|\varepsilon_N^0\|_{\widehat{N}}
   +c_3\kappa N^{-m}\max\limits_{0\leq l\leq k}\|f^l\|_m\\
   &\quad+c_7(N-1)^{1-m}\max\limits_{0\leq l\leq k}\|u^l\|_m,\quad k\geq 2.
\end{align*}
Similar to the proof of Theorem \ref{thm-stable-sem-dis}, we can immediately get the following conclusions:
\begin{eqnarray*}
   && \|\varepsilon_N^{k}\|_{1,{\widehat{N}}} \leq\|\varepsilon_N^0\|_{\widehat{N}}+
   \frac{1}{\zeta_{1,k}}\left(c_3\kappa N^{-m}\max\limits_{0\leq l\leq k}\|f^l\|_m+c_7(N-1)^{1-m}\max\limits_{0\leq l\leq k}\|u^l\|_m\right).
\end{eqnarray*}
Notice that
\begin{eqnarray*}
   && \|\varepsilon_N^0\|_{\widehat{N}} = \|u^0-u_N^0\|_{\widehat{N}}=
  \|u^0(x)-I_Nu^0(x)\|_{\widehat{N}}=0,
\end{eqnarray*}
and the boundedness of $\kappa$ and $\zeta_{1,k}^{-1}$, then there exists a constant $\overline{c}>0$ such that
\begin{eqnarray*}
   && \|\varepsilon_N^{k}\|_{1,{\widehat{N}}} \leq
   \overline{c}\exp\left(\frac{\alpha_nT}{1-\alpha_n}\right)\left( N^{-m}\max\limits_{0\leq l\leq k}\|f^l\|_m+(N-1)^{1-m}\max\limits_{0\leq l\leq k}\|u^l\|_m\right).
\end{eqnarray*}
\hfill$\Box$
\section{Numerical validation} \label{sec:4}
\subsection{Implementation}
We provide a comprehensive account of the implementation of problem \eqref{var-full-dis-k1-coll} using the shifted Legendre collocation method.

Considering problem (\ref{var-full-dis-k1-coll}), we express the function $u_N^{k+1}$ in terms of the Lagrangian interpolants based on the shifted Legendre-Gauss-Lobatto points $\hat{\xi}_{i},~i=0,1,\cdots,N$, i.e.,
\begin{eqnarray}
   &&   u_N^{k}=\sum\limits_{i=0}^Nc_{i}^{k}h_i(x),\label{expr-u-lagr}
\end{eqnarray}
where $c_{i}^{k}:=u_N^{k}(\hat{\xi}_i)$, unknowns of the discrete solution. $h_i(x)$ is the Lagrangian polynomials defined in $\Lambda$, which satisfies
\begin{eqnarray*}
   && h_i(x)\in\mathds{P}_N(\Lambda),\quad h_i(\hat{\xi}_j)=\delta_{ij},\quad i,j=0,1,\cdots,N,
\end{eqnarray*}
where $\delta_{ij}$ is the Kronecker symbols. Taking (\ref{expr-u-lagr}) into (\ref{var-full-dis-k1-coll}), and notice that the homogeneous Dirichlet boundary condition (\ref{pro-bou-cond}), then choosing each test function $v_N$ to be $h_l(x)$~$(l=1,2,\cdots,N-1)$, we have
\begin{eqnarray*}
   && \sum\limits_{i=1}^{N-1}\left(h_i(x),h_l(x)\right)_{\widehat{N}}c_i^k+\kappa\sum\limits_{i=1}^{N-1}\left(\frac{\mathrm{d}h_i(x)}{\mathrm{d}x},\frac{\mathrm{d}h_l(x)}{\mathrm{d}x}\right)_{\widehat{N}}
   c_i^k=F_{\widehat{N}}^k\left(h_l(x)\right),\quad k=1,2,\cdots,N_T.
\end{eqnarray*}
Define the matrices
\begin{eqnarray*}
   && R:=\left(r_{ij}\right)_{i,j=1}^{N-1},\qquad r_{ij}:=\left(h_i(x),h_j(x)\right)_{\widehat{N}}=\sum\limits_{p=0}^{N}h_i(\hat{\xi_p})h_j(\hat{\xi_p})\hat{\omega}_p=\delta_{ij}\hat{\omega}_i,\\
   && G:=\left(g_{ij}\right)_{i,j=1}^{N-1},\qquad g_{ij}:=\left(\frac{\mathrm{d}h_i(x)}{\mathrm{d}x},\frac{\mathrm{d}h_l(x)}{\mathrm{d}x}\right)_{\widehat{N}}=\sum\limits_{p=0}^{N}\frac{\mathrm{d}h_i(\hat{\xi_p})}{\mathrm{d}x}\frac{\mathrm{d}h_j(\hat{\xi_p})}{\mathrm{d}x}\hat{\omega}_p,\\
   && C^{k}:=\left(c_1^k,c_2^k,\cdots,c_{N-1}^k\right)^T,\\
   && Q^{k}:=\left(F_{\widehat{N}}^{k}(h_1(x)),F_{\widehat{N}}^{k}(h_2(x)),\cdots,F_{\widehat{N}}^{k}(h_{N-1}(x))\right)^T.
\end{eqnarray*}
Then, we obtain the matrix representation of the above equation in the following form:
\begin{eqnarray}
   && \left(R+\kappa G\right)C^k=Q^k,\quad k=1,2,\cdots,N_T.\label{coll-matrix}
\end{eqnarray}
The linear system \eqref{coll-matrix} can be solved in particular by the LU factorization or other related computational techniques.

Finally, we discuss about the calculation of $Q^{k}$. When $k=0$, the initial condition $u_N^0$ taken to be
\begin{eqnarray*}
   && u_N^{0}=\sum\limits_{i=0}^Nc_{i}^{0}h_i(x),\quad c_{i}^0=u^0\left(\hat{\xi}_i\right).
\end{eqnarray*}
It is not difficult to see that
\begin{eqnarray*}
   && u_N^0\left(\hat{\xi}_i\right)=u^0\left(\hat{\xi}_i\right),\quad i=0,1,\cdots,N,
\end{eqnarray*}
which implies that $u_N^0$ satisfies interpolation condition (\ref{inter-cond}). When $k\geq 1$, suppose that
\begin{eqnarray*}
   && u_N^m=\sum\limits_{i=1}^{N-1}c_{i}^{m}h_i(x),\quad m=1,2,\cdots,k,
\end{eqnarray*}
then
\begin{eqnarray*}
   && \left(u_N^m,h_l(x)\right)_{\widehat{N}}=\sum\limits_{i=1}^{N-1}c_{i}^{m}\sum\limits_{p=0}^{N}
   h_i(\hat{\xi_p}) h_l(\hat{\xi_p})\hat{\omega}_{p}=c_{l}^m \hat{\omega}_{l},\quad l=1,2,\cdots,N-1.
\end{eqnarray*}
Furthermore,
\begin{eqnarray*}
   &&  \left(f^{k},h_l(x)\right)_{\widehat{N}}=\sum\limits_{p=0}^{N}f^{k}(\hat{\xi_p})h_l(\hat{\xi_p}))\hat{\omega}_{p}
   =f^{k}(\hat{\xi_l})\hat{\omega}_{l},\quad l=1,2,\cdots,N-1.
\end{eqnarray*}
In a word, we can easily obtain $Q^{k}$  at each iteration of time-step.

\subsection{Numerical results}
We present a series of numerical results to validate our theoretical propositions.

Firstly, to investigate the computational performance of two discrete fractional differential operators $L_t^{\alpha}$ and $F_t^{\alpha}$, we test three examples from \cite{akman2018new}. Denote $\beta:=\frac{\alpha}{1-\alpha}$.
\begin{examp}\label{examp-1}
Consider the function $h(t)=t^m$~$(m\geq 1,~m\in\mathds{N})$, the Caputo-Fabrizio fractional derivative of
order $\alpha$ with $0<\alpha<1$ of $h(t)$ is written as
\begin{eqnarray*}
   && {}_0^{\mathrm{CF}}\!D_t^{\alpha}t^m=\frac{1}{1-\alpha}\left\{\sum_{i=0}^{m-1}(-1)^i\frac{m!}{(m-i-1)!\beta^{i+1}}t^{m-i-1}+(-1)^m\exp(-\beta t)\frac{m!}{\beta^m}\right\}.
\end{eqnarray*}
\end{examp}

\begin{examp}\label{examp-2}
Consider the function $h(t)=\cos(\omega t)$, the Caputo-Fabrizio fractional derivative of
order $\alpha$ with $0<\alpha<1$ of $h(t)$ is written as
\begin{eqnarray*}
   &&  {}_0^{\mathrm{CF}}\!D_t^{\alpha}\cos(\omega t)=-\frac{1}{1-\alpha}\frac{\beta^2\omega}{\beta^2+\omega^2}\left\{\frac{\sin(\omega t)}{\beta}-\omega\frac{\cos(\omega t)}{\beta^2}+\exp(-\beta t)\frac{\omega}{\beta^2} \right\}.
\end{eqnarray*}
\end{examp}

\begin{examp}\label{examp-3}
Consider the function $h(t)=\exp(\omega t)$, the Caputo-Fabrizio fractional derivative of
order $\alpha$ with $0<\alpha<1$ of $h(t)$ is written as
\begin{eqnarray*}
   &&  {}_0^{\mathrm{CF}}\!D_t^{\alpha}\exp(\omega t)=\frac{1}{1-\alpha}\omega\frac{\exp(\omega t)-\exp(-\beta t)}{\omega+\beta}.
\end{eqnarray*}
\end{examp}
\noindent The proofs based on the method of integration by parts can be found in \cite{caputo2015new} and \cite{akman2018new}.

We choose $m=4$ in Example \ref{examp-1} and $\omega=5$ in Example \ref{examp-2} and  Example \ref{examp-3}, and set $\alpha=0.5$,~$t\in[0,2]$. Define the errors
\begin{eqnarray*}
   && E_1(\Delta t):=\left|{}_0^{\mathrm{CF}}\!D_t^{\alpha}h^{N_T}-L_t^{\alpha}h^{N_T}\right|,\quad  E_2(\Delta t):=\left|{}_0^{\mathrm{CF}}\!D_t^{\alpha}h^{N_T}-F_t^{\alpha}h^{N_T}\right|,
\end{eqnarray*}
for $L_t^{\alpha}$ and $F_t^{\alpha}$ operators, respectively, where $N_T$ is the last time step. Tables \ref{tab-1}--\ref{tab-3} give the
numerical results of approximation error and CPU time with three examples. Here CPU time represents the total computation time, that is, the whole time for computing the approximations of Caputo-Fabrizio fractional derivatives at every time step. The convergence rates in Tables are given by
\begin{eqnarray*}
   && \mathrm{Rate}^1:=\log_2\left(\frac{E_1(2\Delta t)}{E_1(\Delta t)}\right),\quad \mathrm{Rate}^2:=\log_2\left(\frac{E_2(2\Delta t)}{E_2(\Delta t)}\right),\quad \mathrm{Rate}^c:=\log_2\left(\frac{\mathrm{CPU}(2N_T)}{\mathrm{CPU}(N_T)}\right).
\end{eqnarray*}
Tables \ref{tab-1}--\ref{tab-3} demonstrate that the errors of both $L_t^{\alpha}$ and $F_t^{\alpha}$ approximations are virtually identical, as a result of their equivalence ($F_t^{\alpha}h^k=L_t^{\alpha}h^k$).
Moreover, both approximations have achieved second-order convergence of error, as stated in Lemma \ref{lem-erro-Ft}.
However, we observe that the CPU time of $F_t^{\alpha}$ approximation increases linearly with respect to $N_T$, while the $L_t^{\alpha}$ approximation increases almost quadratically.
This suggests that the $F_t^{\alpha}$ operator holds promise as it requires less storage and incurs lower computational costs than the $L_t^{\alpha}$ operator during computation.

\begin{table}[!h]\small
\centering
\caption{Comparisons of $L_t^{\alpha}$ with $F_t^{\alpha}$ for Example \ref{examp-1}.}\label{tab-1}
\begin{tabular}{lllllllll}
  \hline
$\Delta t$ & $E_1(\Delta t)$ & $\mathrm{Rate}^1$ & $\mathrm{CPU}$ & $\mathrm{Rate}^c$ & $E_2(\Delta t)$ & $\mathrm{Rate}^2$ & $\mathrm{CPU}$ & $\mathrm{Rate}^c$   \\\hline

$2/5000$ & $5.5339e-7$ & $-$ & $0.2863$ & $-$ & $5.5339e-7$ & $-$ & $0.0137$ & $-$\\

$2/10000$ & $1.3835e-7$ & $2.0000$ & $1.1262$ & $1.9759$ & $1.3835e-7$ & $2.0000$ & $0.0242$ & $0.8202$\\

$2/20000$ & $3.4586e-8$ & $1.9996$ & $4.4039$ & $1.9673$ & $3.4586e-8$ & $1.9996$ & $0.0456$ & $0.9140$\\

$2/40000$ & $8.6339e-9$ & $2.0025$ & $17.6499$ & $2.0028$ & $8.6340e-9$ & $2.0025$ & $0.0832$ & $0.8677$\\

$2/80000$ & $2.1765e-9$ & $1.9880$ & $72.2386$ & $2.0331$ & $2.1764e-9$ & $1.9881$ & $0.1675$ & $1.0100$\\
\hline
\end{tabular}
\end{table}

\begin{table}[!h]\small
\centering
\caption{Comparisons of $L_t^{\alpha}$ with $F_t^{\alpha}$ for Example \ref{examp-2}.}\label{tab-2}
\begin{tabular}{lllllllll}
  \hline
$\Delta t$ & $E_1(\Delta t)$ & $\mathrm{Rate}^1$ & $\mathrm{CPU}$ & $\mathrm{Rate}^c$ & $E_2(\Delta t)$ & $\mathrm{Rate}^2$ & $\mathrm{CPU}$ & $\mathrm{Rate}^c$   \\\hline

$2/5000$ & $9.4713e-8$ & $-$ & $0.2928$ & $-$ & $9.4713e-8$ & $-$ & $0.0118$ & $-$\\

$2/10000$ & $2.3683e-8$ & $2.0000$ & $1.1333$ & $1.9525$ & $2.3683e-8$ & $2.0000$ & $0.0236$ & $1.0000$\\

$2/20000$ & $5.9207e-9$ & $2.0000$ & $4.3878$ & $1.9530$ & $5.9207e-9$ & $2.0000$ & $0.0417$ & $0.8212$\\

$2/40000$ & $1.4808e-9$ & $1.9994$ & $17.6614$ & $2.0090$ & $1.4808e-9$ & $1.9994$ & $0.0912$ & $1.1290$\\

$2/80000$ & $3.6928e-10$ & $2.0036$ & $72.0036$ & $2.0275$ & $3.6933e-10$ & $2.0034$ & $0.1597$ & $0.8083$\\
\hline
\end{tabular}
\end{table}

\begin{table}[!h]\small
\centering
\caption{Comparisons of $L_t^{\alpha}$ with $F_t^{\alpha}$ for Example \ref{examp-3}.}\label{tab-3}
\begin{tabular}{lllllllll}
  \hline
$\Delta t$ & $E_1(\Delta t)$ & $\mathrm{Rate}^1$ & $\mathrm{CPU}$ & $\mathrm{Rate}^c$ & $E_2(\Delta t)$ & $\mathrm{Rate}^2$ & $\mathrm{CPU}$ & $\mathrm{Rate}^c$   \\\hline

$2/5000$ & $2.4474e-3$ & $-$ & $0.2927$ & $-$ & $2.4474e-3$ & $-$ & $0.0171$ & $-$\\

$2/10000$ & $6.1184e-4$ & $2.0000$ & $1.1284$ & $1.9468$ & $6.1184e-4$ & $2.0000$ & $0.0325$ & $0.9264$\\

$2/20000$ & $1.5296e-4$ & $2.0000$ & $4.4847$ & $1.9907$ & $1.5296e-4$ & $2.0000$ & $0.0593$ & $0.8676$\\

$2/40000$ & $3.8215e-5$ & $2.0001$ & $18.0631$ & $2.0100$ & $3.8215e-5$ & $2.0001$ & $0.1034$ & $0.8021$\\

$2/80000$ & $9.5891e-6$ & $1.9947$ & $74.0835$ & $2.0231$ & $9.5891e-6$ & $1.9947$ & $0.1700$ & $0.7173$\\
\hline
\end{tabular}
\end{table}

Secondly, we provide preliminary computational findings to demonstrate the efficacy of the finite difference/shifted Legendre collocation method (abbreviated as FCM).

\begin{examp}\label{examp-4}
 Consider the following three-term time-fractional diffusion equations:
\begin{equation}\label{examp1}
   \left\{\begin{array}{ll}
            \sum\limits_{i=1}^3d_i\cdot{}_0^{\mathrm{CF}}\!D_t^{\alpha_{i}}u(x,t)=\frac{\partial^2 u(x,t)}{\partial x^2}+f(x,t), & (x,t)\in\Lambda\times(0,1], \\
            u(x,0)=\sin x, & x\in\Lambda, \\
            u(0,t)=u(\pi,t)=0, & t\in[0,1],
          \end{array}\right.
\end{equation}
where $\Lambda=(0,\pi)$, $0<\alpha_1\leq\alpha_2\leq\alpha_3<1$, and
 \begin{eqnarray*}
    && f(x,t)=\sin x \left\{ 1+t^2+\sum\limits_{i=1}^3\frac{d_i}{1-\alpha_i}\left[ \frac{2}{\beta_i}t-\frac{2}{\beta_i^2}+\exp(-\beta_it)\frac{2}{\beta_i^2} \right] \right\},
 \end{eqnarray*}
with $\beta_i=\frac{\alpha_i}{1-\alpha_i}$,~$i=1,2,3$. The exact solution of the Eq.\eqref{examp1} is $u(x,t)=(1+t^2)\sin x$, which is sufficiently smooth. In our experiments, we set the parameters $d_1=1,d_2=2,d_3=3$.
\end{examp}

The following error norms have been used as the error indicator:
\begin{eqnarray*}
   && \|e\|_\infty:=\|u^{k}(x)-u_N^{k}\|_{L^\infty(\Lambda)}=\sup\limits_{x\in\Lambda}|u^{k}(x)-u_N^{k}|,\\
   && \|e\|_0:=\|u^{k}(x)-u_N^{k}\|_0,\\
   &&\|e\|_1:=\|u^{k}(x)-u_N^{k}\|_1.
\end{eqnarray*}
We test Example \ref{examp-4} with four cases: Case 1. $\alpha_1=\alpha_2=\alpha_3=0.5$, then (\ref{examp1}) reduces to the single-term time-fractional diffusion equation; Case 2. $\alpha_1=0.3$,~$\alpha_2=0.5$,~$\alpha_3=0.7$; Case 3. $\alpha_1=0.2$,~$\alpha_2=0.3$,~$\alpha_3=0.4$; Case 4. $\alpha_1=0.6$,~$\alpha_2=0.7$,~$\alpha_3=0.8$. In time discretization, we use $F_t^{\alpha_i}$ operator. Table \ref{tab-4} shows the errors and temporal accuracy of FCM with polynomial degree $N=20$ at $t_k=1$ for different cases of $\alpha_i$. Here convergence rates are given by
\begin{eqnarray*}
   && \mathrm{Rate}=\log_2\left(\frac{\|e(N,2\Delta t)\|_l}{\|e(N,\Delta t)\|_l}\right),\qquad l=\infty,0,1.
\end{eqnarray*}
It can be observed that the FCM  exhibits a second-order temporal convergence rate, , which is consistent with our theoretical analysis.

\begin{table}[!h]\small
\centering
\caption{Numerical convergence of FCM in the temporal direction for Example \ref{examp-4}.}\label{tab-4}
\begin{tabular}{llllllll}
  \hline
$\alpha_i$ & $\Delta t$ & $\|e\|_\infty$ & $\mathrm{Rate}$ & $\|e\|_0$ & $\mathrm{Rate}$ & $\|e\|_1$ & $\mathrm{Rate}$   \\\hline

 & $1/160$ & $2.7461e-5$ & $-$ & $3.4423e-5$ & $-$ & $4.8685e-5$ & $-$  \\

$\alpha_1=0.5$ & $1/320$ & $6.8824e-6$ & $1.9964$ & $8.6272e-6$ & $1.9964$ & $1.2201e-5$ & $1.9965$  \\

 $\alpha_2=0.5$ & $1/640$ & $1.7227e-6$ & $1.9982$ & $2.1595e-6$ & $1.9982$ & $3.0541e-6$ & $1.9982$  \\

$\alpha_3=0.5$ & $1/1280$ & $4.3095e-7$ & $1.9991$ & $5.4020e-7$ & $1.9991$ & $7.6400e-7$ & $1.9991$ \\

 & $1/2560$ & $1.0777e-7$ & $1.9996$ & $1.3509e-7$ & $1.9996$ & $1.9106e-7$ & $1.9995$  \\
\hline

 & $1/160$ & $2.3560e-5$ & $-$ & $2.9532e-5$ & $-$ & $4.1767e-5$ & $-$  \\

$\alpha_1=0.3$ & $1/320$ & $5.9153e-6$ & $1.9938$ & $7.4149e-6$ & $1.9938$ & $1.0487e-5$ & $1.9938$  \\

 $\alpha_2=0.5$ & $1/640$ & $1.4820e-6$ & $1.9969$ & $1.8577e-6$ & $1.9969$ & $2.6274e-6$ & $1.9969$  \\

$\alpha_3=0.7$ & $1/1280$ & $3.7090e-7$ & $1.9984$ & $4.6493e-7$ & $1.9984$ & $6.5755e-7$ & $1.9985$ \\

 & $1/2560$ & $9.2776e-8$ & $1.9992$ & $1.1630e-7$ & $1.9991$ & $1.6448e-7$ & $1.9992$  \\
\hline

 & $1/160$ & $3.0350e-5$ & $-$ & $3.8047e-5$ & $-$ & $5.3810e-5$ & $-$  \\

$\alpha_1=0.2$ & $1/320$ & $7.5961e-6$ & $1.9984$ & $9.5218e-6$ & $1.9985$ & $1.3467e-5$ & $1.9984$  \\

 $\alpha_2=0.3$ & $1/640$ & $1.9000e-6$ & $1.9993$ & $2.3817e-6$ & $1.9992$ & $3.3684e-6$ & $1.9993$  \\

$\alpha_3=0.4$ & $1/1280$ & $4.7513e-7$ & $1.9996$ & $5.9558e-7$ & $1.9996$ & $8.4233e-7$ & $1.9996$ \\

 & $1/2560$ & $1.1880e-7$ & $1.9998$ & $1.4892e-7$ & $1.9998$ & $2.1061e-7$ & $1.9998$  \\
\hline

 & $1/160$ & $1.4817e-5$ & $-$ & $1.8573e-5$ & $-$ & $2.6268e-5$ & $-$  \\

$\alpha_1=0.6$ & $1/320$ & $3.7588e-6$ & $1.9789$ & $4.7117e-6$ & $1.9789$ & $6.6637e-6$ & $1.9789$  \\

 $\alpha_2=0.7$ & $1/640$ & $9.4656e-7$ & $1.9895$ & $1.1865e-6$ & $1.9895$ & $1.6781e-6$ & $1.9895$  \\

$\alpha_3=0.8$ & $1/1280$ & $2.3750e-7$ & $1.9948$ & $2.9771e-7$ & $1.9947$ & $4.2105e-7$ & $1.9948$ \\

 & $1/2560$ & $5.9483e-8$ & $1.9974$ & $7.4563e-8$ & $1.9974$ & $1.0545e-7$ & $1.9974$  \\
\hline
\end{tabular}
\end{table}

\begin{figure}[!h]
  \centering\caption{Numerical convergence of FCM in the spatial direction for Example \ref{examp-4}.}\label{Fig1}
 \centerline{\scalebox{0.25}{\includegraphics{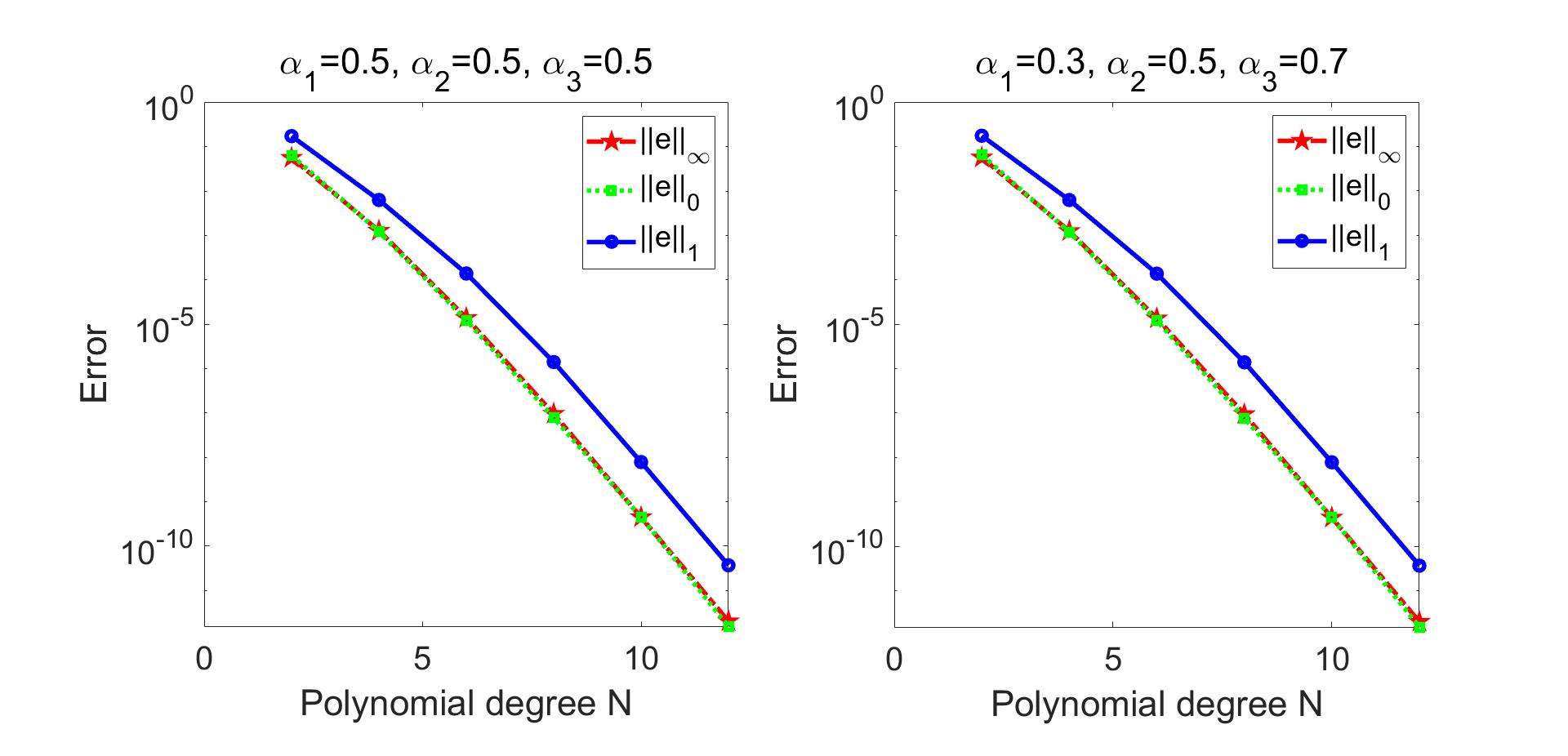}}}\vspace{-10pt}
  \centerline{\scalebox{0.25}{\includegraphics{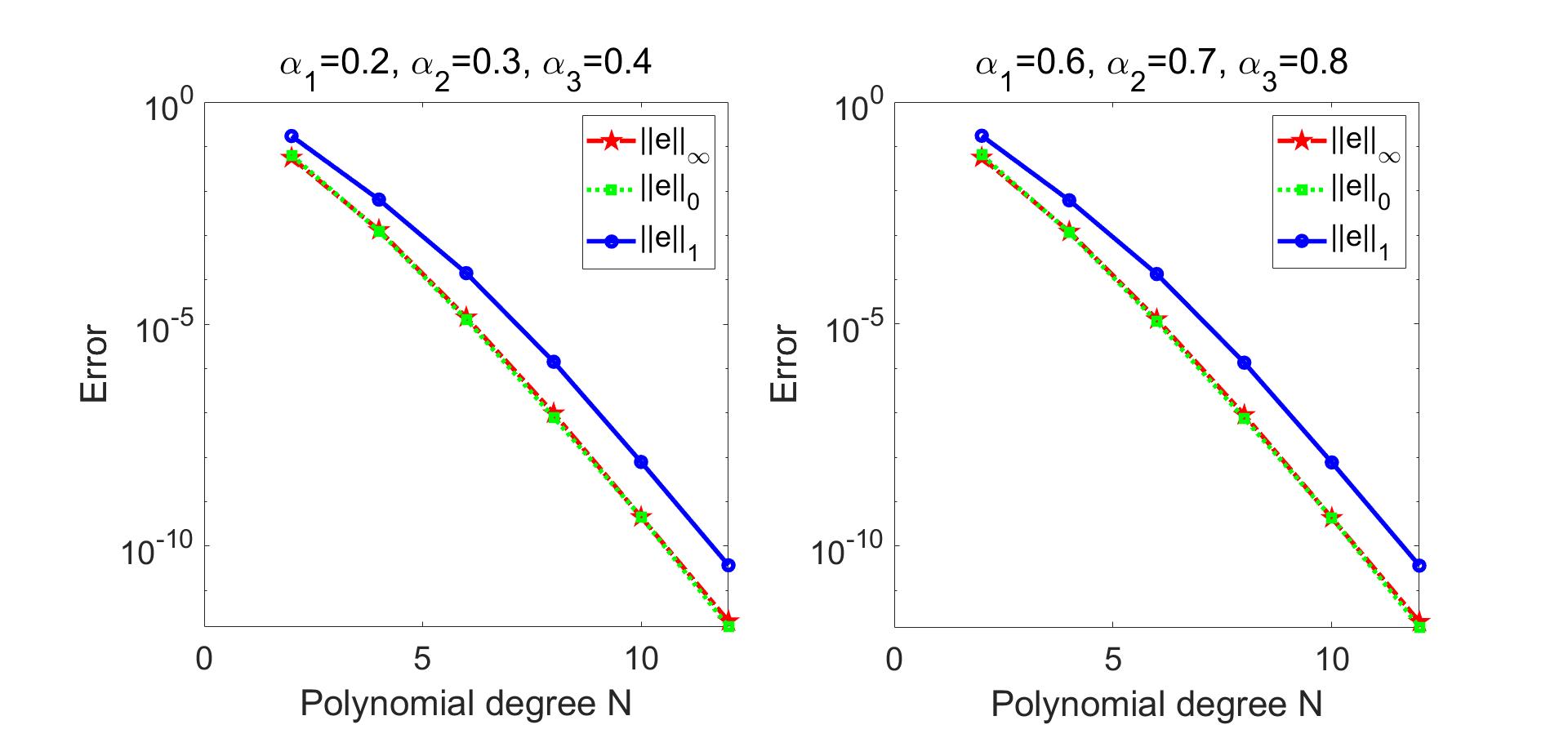}}}
\end{figure}

Next, we check the spatial accuracy with respect to the polynomial degree $N$. In order to avoid the contamination of temporal error, we need
fix the time step $\Delta t$ sufficiently small. Here we take $\Delta t=10^{-6}$, and terminate computing at $t_k=0.01$ for saving time.
Fig. \ref{Fig1} shows the errors with respect to polynomial degree $N$ at $t_k=0.01$ in semi-log scale. Evidently, the spatial discretization exhibits exponential convergence as demonstrated by the nearly linear curves depicted in this figure. The aforementioned is known as spectral accuracy as expected since the exact solution  is a sufficiently smooth function with respect to the space variable.

To further verify the numerical validity, we finally  test a two-dimensional problem.
\begin{examp}\label{examp-5}
 Consider the following three-term time-fractional diffusion equations:
\begin{equation}\label{examp5}
   \left\{\begin{array}{ll}
            \sum\limits_{i=1}^3d_i\cdot{}_0^C\!D_t^{\alpha_{i}}u(x,y,t)=\Delta u(x,y,t)+f(x,y,t), & (x,y,t)\in\Omega\times(0,1], \\
            u(x,y,0)=\sin (2\pi xy), & (x,y)\in\Omega, \\
            u(x,y,t)=0, & (x,y,t)\in\partial\Omega\times(0,1],
          \end{array}\right.
\end{equation}
where $\Omega=[0,1]\times[0,1]$, $0<\alpha_1\leq\alpha_2\leq\alpha_3\leq1$, and
 \begin{eqnarray*}
    && f(x,y,t)=4\pi^2(x^2+y^2)\sin(2\pi xy)\exp(t)+\sum\limits_{i=1}^3\frac{d_i}{1-\alpha_i}\frac{\exp(t)-\exp(\beta_it)}{1+\beta}.
 \end{eqnarray*}
The exact solution of the Eq. (\ref{examp5}) is $u(x,y,t)=\sin(2\pi xy)\exp(t)$, which is sufficiently smooth. In our experiments, we set the parameters $d_1=1,d_2=2,d_3=3$.
\end{examp}

\begin{figure}[!h]
  \centering\caption{Numerical convergence of FCM in the spatial direction for Example \ref{examp-5}.}\label{Fig2}
 \centerline{\scalebox{0.25}{\includegraphics{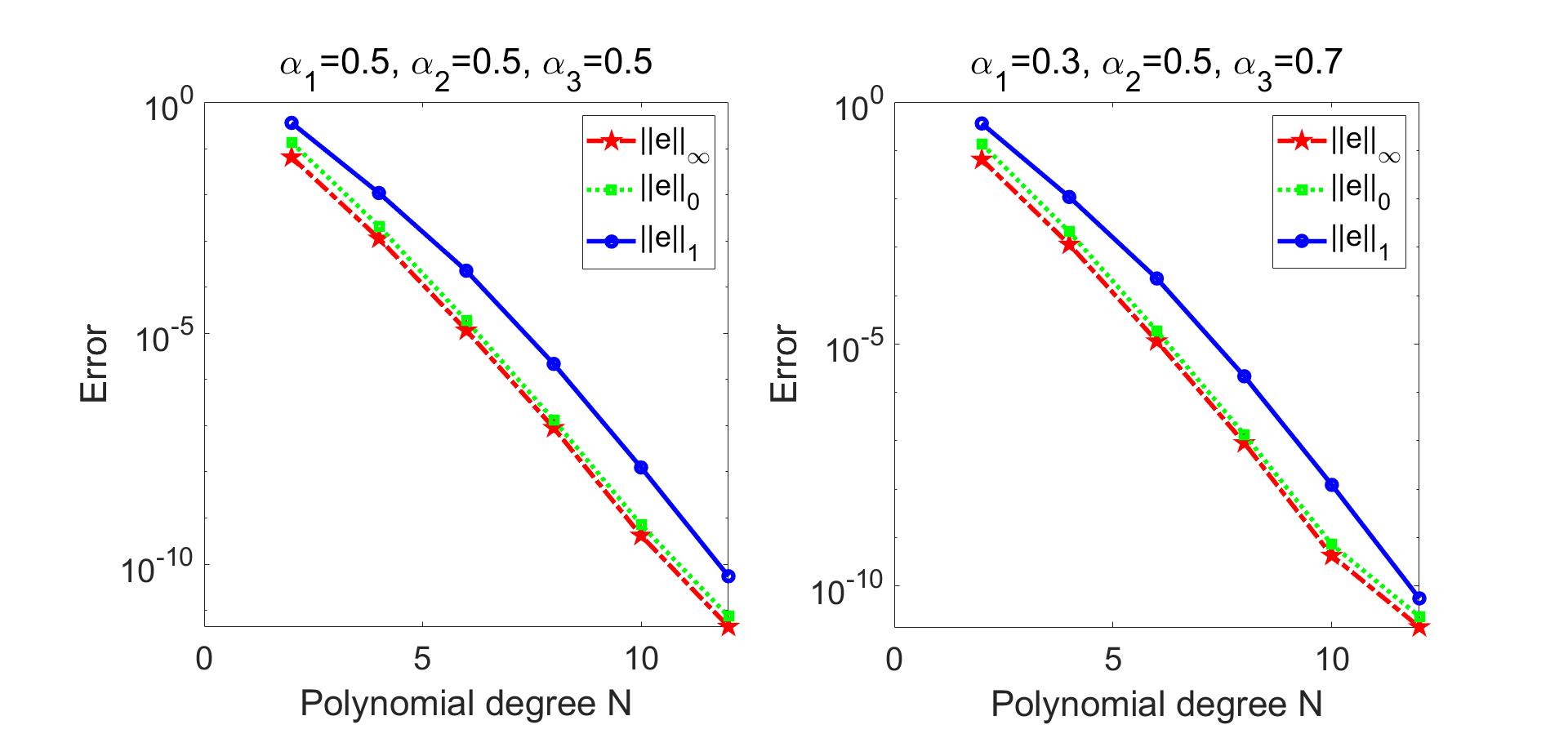}}}\vspace{-6pt}
  \centerline{\scalebox{0.25}{\includegraphics{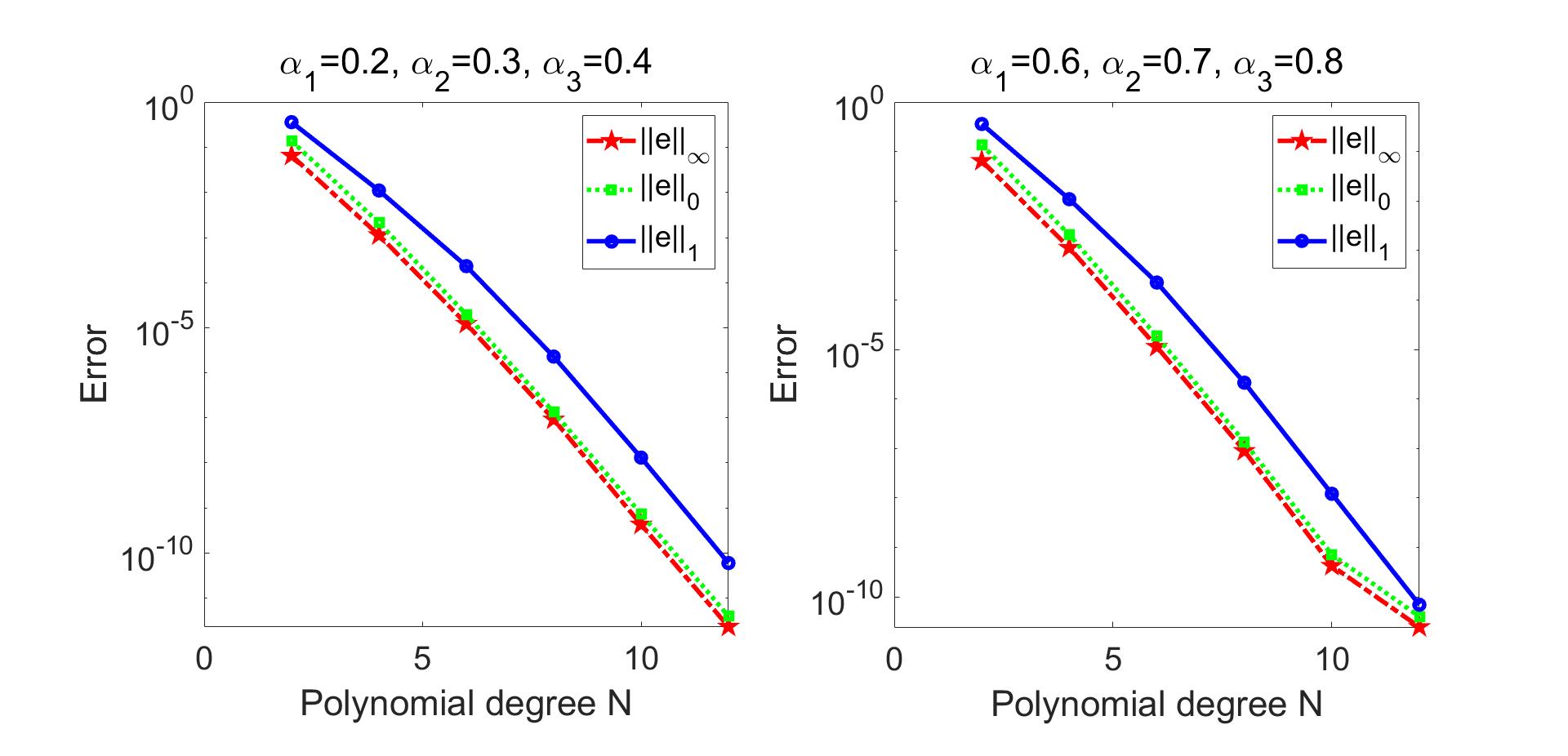}}}
\end{figure}

In this case, we denote $\{\hat{\xi}_{pq}\}_{p,q=0}^N:=\{(\hat{\xi}_p,\hat{\xi}_q)\}_{p,q=0}^N$ and $\{\hat{\omega}_{pq}\}_{p,q=0}^N:=\{\hat{\omega}_p\hat{\omega}_q\}_{p,q=0}^N$ as the nodes and weights of shifted Legendre-Gauss-Lobatto quadratures on $\overline{\Omega}$. Then we express the function $u_N^{k+1}$ in terms of the two-dimensional Lagrangian interpolants based on the shifted Legendre-Gauss-Lobatto points $\hat{\xi}_{ij},~i,j=0,1,\cdots,N$,
\begin{equation}\label{expr-u-lagr}
    u_N^{k+1}=\sum\limits_{i=0}^N\sum\limits_{j=0}^Nc_{ij}^{k+1}h_i(x)h_j(y),
\end{equation}
where $c_{ij}^{k+1}:=u_N^{k+1}(\hat{\xi}_i,\hat{\xi}_j)$, unknowns of the discrete solution. $h_i(x)$ and $h_j(y)$ are the Lagrangian polynomials defined in $I_x:=[0,1]$ and $I_y:=[0,1]$, i.e.,
\begin{eqnarray*}
   && h_i(x)\in\mathrm{P}_N(I_x),\quad h_i(\hat{\xi}_l)=\delta_{il},\quad i,l=0,1,\cdots,N,\\
   && h_j(y)\in\mathrm{P}_N(I_y),\quad h_j(\hat{\xi}_s)=\delta_{js},\quad j,s=0,1,\cdots,N,
\end{eqnarray*}
where $\delta_{il}$ and $\delta_{js}$ are the Kronecker symbols. A linear system such as \eqref{coll-matrix} can be readily derived.
Here we take $\Delta t=10^{-6}$. Fig. \ref{Fig2} shows the errors with respect to polynomial degree $N$ in semi-log scale.
Thanks to the fast scheme \eqref{app-CFD-rec-tk}, a small time step does not significantly escalate the computational burden in the time direction,
thereby the proposed method is effective even for handling high-dimensional problems.

\section{Concluding remarks}\label{sec:5}

In this work, we have developed a fully discrete scheme for the multi-term time-fractional diffusion
equations with Caputo-Fabrizio derivatives. The proposed approach utilizes the finite difference method to approximate multi-term fractional derivatives in time and employs the Legendre spectral collocation method for spatial discretization.
Specifically, we use the exponential property of Caputo-Fabrizio derivative to give a recursive difference calculation scheme, which offers benefits in terms of computational complexity and storage capacity.
The proposed scheme has been proved to be unconditionally stable and convergent with order $O\left(\left(\Delta t\right)^{2}+N^{-m}\right)$. Numerical results show good agreement with the theoretical analysis. Due to its high resolution feature in spectral approximation, the proposed method can be extended to handle multi-term time-fractional diffusion equations in higher spatial dimensions.

%%%%%%%%%%%%%%%%%%%%%%%%%%%%%%%%%%%%%%%%%%%%%%%%%%%%%%
%          AI TOOLS, USE AND LOCATION
%%%%%%%%%%%%%%%%%%%%%%%%%%%%%%%%%%%%%%%%%%%%%%%%%%%%%%
%We follow COPE's guidelines and policies regarding the use of Artificial Intelligence (AI) tools. COPE Policy on AI tools can be found at https://publicationethics.org/cope-position-statements/ai-author.

%Authors using AI tools in the writing of a manuscript, production of images or graphical elements of the paper, or in the collection and analysis of data, must be transparent in disclosing in this section how the AI tool was used and which tool was used. Authors are fully responsible for the content of their manuscript, even those parts produced by an AI tool, and are thus liable for any breach of publication ethics. - COPE

%Disclosure instructions

%If there is nothing to disclose, there is no need to add a declaration, otherwise please declare.

%\section*{Use of AI tools declaration}
%The author(s) declare(s) they have used Artificial Intelligence (AI) tools in the creation of this article.
%AI tools used:
%How were the AI tools used?
%Where in the article is the information located?

\section*{Conflict of interest}
The authors declare no conflicts of interest in this manuscript.

%\bibliographystyle{plain}

%\bibliography{ref_fastCF_2023}

\end{document}